\documentclass{gen-j-l}
\usepackage{amssymb}
\usepackage{amsfonts}

\newtheorem{theorem}{Theorem}[section]
\newtheorem{lemma}[theorem]{Lemma}
\theoremstyle{definition}
\newtheorem{definition}[theorem]{Definition}
\newtheorem{example}[theorem]{Example}

\newtheorem{remark}[theorem]{Remark}
\numberwithin{equation}{section}
\theoremstyle{plain}

\newtheorem{corollary}[theorem]{Corollary}

\newtheorem{proposition}[theorem]{Proposition}

\copyrightinfo{2001}{enter name of copyright holder}
\input{tcilatex}

\begin{document}
\title[Zeta functions and Newton polyhedra]{Zeta functions for analytic
mappings, log-principalization of ideals, and Newton polyhedra}
\author{Willem Veys}
\address{University of Leuven, Department of Mathematics\\
Celestijnenlaan 200 B, B-3001 Leuven (Heverlee), Belgium}
\email{wim.veys@wis.kuleuven.be}
\thanks{The first author was partially supported by the Fund of Scientific
Research - Flanders (G.0318.06).}
\author{W. A. Zuniga-Galindo}
\address{Department of Mathematics and Computer Science, Barry University,
11300 N.E. Second Avenue, Miami Shores, Florida 33161, USA}
\email{wzuniga@mail.barry.edu}
\thanks{The second author thanks the financial support of the NSA. Project
sponsored by the National Security Agency under Grant Number
H98230-06-1-0040. The United States Government is authorized to reproduce
and distribute reprints notwithstanding any copyright notation herein.}

\begin{abstract}
In this paper we provide a geometric description of the possible poles of
the Igusa local zeta function $Z_{\Phi }(s,\mathbf{f})$ associated to an
analytic mapping $\mathbf{f}=$ $\left( f_{1},\ldots ,f_{l}\right)
:U(\subseteq K^{n})\rightarrow K^{l}$, and a locally constant function $\Phi 
$, with support in $U$, in terms of a log-principalizaton of the $K\left[ x%
\right] -$ideal $\mathcal{I}_{\mathbf{f}}=\left( f_{1},\ldots ,f_{l}\right) $%
. Typically our new method provides a much shorter list of possible poles
compared with the previous methods. We determine the largest real part of
the poles of the Igusa zeta function, and then as a corollary, we \ obtain \
an asymptotic estimation for the number of solutions of an arbitrary system
of polynomial congruences in terms of the log-canonical threshold \ of the
subscheme \ given by $\mathcal{I}_{\mathbf{f}}$. We associate to an analytic
mapping $\boldsymbol{f}$ $=$ $\left( f_{1},\ldots ,f_{l}\right) $ a Newton
polyhedron $\Gamma \left( \boldsymbol{f}\right) $ and \ a new notion of
non-degeneracy with respect to $\Gamma \left( \boldsymbol{f}\right) $. The
novelty of this notion resides in the fact that it depends on one Newton
polyhedron, and Khovanskii's non-degeneracy notion depends on the Newton
polyhedra of $f_{1},\ldots ,f_{l}$ . By \ constructing a
log-principalization, we give an explicit list \ for the possible poles of $%
Z_{\Phi }(s,\mathbf{f})$, $l\geq 1$, \ in the case in which $\mathbf{f}$ \
is non-degenerate with respect to $\Gamma \left( \boldsymbol{f}\right) $.
\end{abstract}

\maketitle

\section{Introduction}

Let $K$ be a $p-$adic field, i.e. $[K:\mathbb{Q}_{p}]<\infty $. Let $R_{K}$\
be the valuation ring of $K$, $P_{K}$ the maximal ideal of $R_{K}$, and $%
\overline{K}=R_{K}/P_{K}$ \ the residue field of $K$. \ The cardinality of
the residue field of $K$ is denoted by $q$, thus $\overline{K}=\mathbb{F}%
_{q} $. For $z\in K$, $ord\left( z\right) \in \mathbb{Z}\cup \{+\infty \}$ \
denotes the valuation of $z$, and $\left\vert z\right\vert
_{K}=q^{-ord\left( z\right) }$ its absolute value. The \ absolute value $%
\left\vert \cdot \right\vert _{K}$ can be extended to $K^{l}$ by defining $%
\left\Vert z\right\Vert _{K}=\max_{1\leq i\leq l}\left\vert z_{i}\right\vert
_{K}$, for $z=(z_{1},\ldots ,z_{l})\in K^{l}$.

Let $f_{1},\ldots ,f_{l}$ be polynomials in $K\left[ x_{1},\ldots ,x_{n}%
\right] $, or, more generally, $K-$analytic functions on an open set $%
U\subset K^{n}$. We consider the mapping $\boldsymbol{f}=$ $\left(
f_{1},\ldots ,f_{l}\right) :K^{n}\rightarrow K^{l}$, respectively, $%
U\rightarrow K^{l}$. Let $\Phi :K^{n}\rightarrow \mathbb{C}$ \ be a
Schwartz-Bruhat function (with support in $U$ in the second case). The Igusa
local zeta function associated to the above data is defined as 
\begin{equation*}
Z_{\Phi }(s,\boldsymbol{f})=Z_{\Phi }(s,\boldsymbol{f},K)=\int%
\limits_{K^{n}}\Phi \left( x\right) \left\Vert \boldsymbol{f}(x)\right\Vert
_{K}^{s}\mid dx\mid ,
\end{equation*}%
for $s\in \mathbb{C}$ with $\func{Re}(s)>0$, where $\mid dx\mid $ is the
Haar measure on $K^{n}$ normalized in such a way that $R_{K}^{n}$ has
measure $1$. \ We write $Z(s,\boldsymbol{f})$, $Z_{0}(s,\boldsymbol{f})$ and 
$Z_{W}(s,\boldsymbol{f})$ when $\Phi $ is the characteristic function of $%
R_{K}^{n}$, $P_{K}^{n}$, and \ an open compact subset $W$ of $K^{n}$,
respectively.

The \ function $Z_{\Phi }(s,\boldsymbol{f})$ admits a meromorphic
continuation to the complex plane as a rational function of $q^{-s}$. Igusa
established this result in the hypersurface case using Hironaka's resolution
theorem \cite[Theorem 8.2.1]{I2}. In the case $l\geq 1$ the rationality of $%
Z_{\Phi }(s,\boldsymbol{f})$ was established by Meuser \ in \cite{M2},
however, as mentioned in the review MR 83g:12015 of \cite{M2}, a trick by
Serre allows to deduce \ the general case from the hypersurface case. Denef
gave a \ completely different proof of the rationality of $Z_{\Phi }(s,%
\boldsymbol{f})$, $l\geq 1$, \ using $p$-adic cell decomposition \cite{D1}.
The mentioned results \ do not give any information about the poles of $%
Z_{\Phi }(s,\boldsymbol{f})$ in the case $l>1.$ In \cite{Z2} the second
author showed that a list of possible poles of \ $Z_{\Phi }(s,\boldsymbol{f}%
) $, $l\geq 1$, can be computed from an embedded resolution of singularities
of the divisor $\cup _{i=1}^{l}f_{i}^{-1}(0)$ by using toroidal geometry. In
the special case in which $\boldsymbol{f}$ is a non-degenerate homogeneous
polynomial mapping the possible poles of $Z_{\Phi }(s,\boldsymbol{f})$ are
given in \cite{Z3}.

In this paper we provide a geometric description of the possible poles of $%
Z_{\Phi }(s,\boldsymbol{f})$, $l\geq 1$, in terms of a log-principalization
of the $K\left[ x\right] -$ideal $\mathcal{I}_{\boldsymbol{f}}=\left(
f_{1},\ldots ,f_{l}\right) $ (see Theorem \ref{Theorem1}). At this point it
is important to mention that the main result in \cite{Z2} gives an algorithm
to compute a list of possible poles of \ $Z_{\Phi }(s,\boldsymbol{f})$, $%
l\geq 1$, in terms of an embedded resolution of singularities of the divisor 
$\cup _{i=1}^{l}f_{i}^{-1}(0)$, while Theorem \ref{Theorem1} gives a list of
candidates to poles in terms of a log-principalization of the ideal $%
\mathcal{I}_{\boldsymbol{f}}$. Typically our new method provides a much
shorter list of possible poles (see Example \ref{Ex1}). It is important to
mention\ that in the case $l=1$ the problem of determining the poles of the
meromorphic continuation of $Z_{\Phi }(s,\boldsymbol{f})$ in $\func{Re}(s)<0$
has been studied extensively (see e.g. \cite{D0}, \cite{I1a}, \cite{S}, \cite%
{M1}, \cite{V1}, \cite{V2}). The relevance of this problem is due to the
existence of several conjectures relating the poles of $Z_{\Phi }(s,%
\boldsymbol{f})$ with \ the structure of the singular locus of $\boldsymbol{f%
}$. In the case of polynomials in two variables, as a consequence of the
works of Igusa, Strauss, Meuser and the first author, there is a complete
solution of this problem \cite{I1a}, \cite{St}, \cite{M1}, \cite{V1a}. For
general polynomials the problem of determination of the poles of $\ Z_{\Phi
}(s,\boldsymbol{f})$ is still open. There exists a generic class of
polynomials named non-degenerate with respect to its Newton polyhedron for
which it is possible to give a small set of candidates for the poles of $%
Z_{\Phi }(s,\boldsymbol{f})$. The poles of the local zeta functions attached
to non-degenerate polynomials can be described in terms of Newton polyhedra.
The case of two variables was studied by Lichtin and Meuser \cite{L-M}. In 
\cite{D1a}, Denef gave a procedure based on monomial changes of variables to
determine a small set of candidates for the poles of $Z_{\Phi }(s,%
\boldsymbol{f})$ in terms of the Newton polyhedron of $f$. This result was
obtained by the second author, using an approach based on the $p-$adic
stationary phase formula and N\'{e}ron $p-$desingularization, for
polynomials with coefficients in a non-archimedean local field of arbitrary
characteristic \cite{Z1}, (see also \cite{D-H}, \cite{S-Z}).

In the case $l=1$, among the conjectures relating the poles of Igusa's zeta
function with topology and singularity theory, we mention here a conjecture
of Igusa that proposes that the real parts of the poles of the Igusa zeta
function of $\boldsymbol{f}$ are roots of the Bernstein polynomial of $%
\boldsymbol{f}$ (see e.g. \cite{D0}, \cite{I2}, and references therein). It
seems reasonable to believe that such relations between poles\ and
singularity theory extend to the case $l>1$. Indeed, recently it was proved
that the\ above-mentioned conjecture of Igusa is valid in the case in which $%
\mathcal{I}_{\boldsymbol{f}}$ is a monomial ideal \cite{H-M-Y}.

In the case $l=1$, the largest real part of the poles of the Igusa zeta\
function has been extensively studied both in the archimedean and
non-archimedean cases \cite{D-H}, \cite{L-M}, \cite{Var}, \cite{Z1}. In the
case $l\geq 1$ we show that the largest real part $-\lambda \left( \mathcal{I%
}_{\boldsymbol{f}}\right) $ of the poles of the Igusa zeta function attached
to $\boldsymbol{f}$ can easily be determined from a log-principalization of
the ideal $\mathcal{I}_{\boldsymbol{f}}$ (see Theorem \ref{theorem1c}). As a
consequence of this result we \ obtain \ an asymptotic estimation for the
number of solutions of an arbitrary system of polynomial congruences in
terms of the log-canonical threshold\ of a log-principalization (see
Corollary \ref{coro1}, and the comments that follow). At this point we have
to mention that in the case $l=1$ Loeser found lower and upper bounds for $%
\lambda \left( \mathcal{I}_{\boldsymbol{f}}\right) $ in terms of certain
geometric invariants introduced by Teissier \cite[Theorem 2.6 and
Proposition 3.1.1]{Lo}, \cite{Te}. In this form he derived a geometric bound
for the number of solutions of a polynomial congruence involving one
polynomial.

If $\boldsymbol{f}$ \ is a polynomial mapping \ with coefficients in a
number field $F$, then for every maximal ideal $P$ of the ring of algebraic
integers of $F$, we can consider $Z(s,\boldsymbol{f},K)$, $l\geq 1$, where $%
K $ is the completion of $F$ with respect to $P$. We give an explicit
formula for $Z(s,\boldsymbol{f},K)$, $l\geq 1$, that is valid for almost all 
$P$ (see Theorem \ref{Theorem1a}). The proof of this formula follows by
adapting the argument given by Denef for the case $l=1$ \cite{D2}.\ 

One can also associate to a sheaf of ideals $\mathcal{I}$ on a smooth
algebraic variety (over a field of characteristic zero) a motivic zeta
function (see \ Definition \ref{defmotzeta}). By using a
log-principalization of $\mathcal{I}$ we give a similar explicit formula for
it\ (see Theorem \ref{Theorem2b}). The proof is a reasonably straightforward
\ generalization of the one given by Denef and Loeser in \ \cite{D-L}. By
specializing to Euler characteristics one obtains the topological zeta
function associated to $\mathcal{I}$.

We attach to an analytic mapping $\boldsymbol{f}$ $=$ $\left( f_{1},\ldots
,f_{l}\right) $ a Newton polyhedron $\Gamma \left( \boldsymbol{f}\right) $
and \ a new notion of non-degeneracy with respect to $\Gamma \left( 
\boldsymbol{f}\right) $. The novelty of this notion resides in the fact that
it depends on \textit{one }Newton polyhedron, and Khovanskii's
non-degeneracy notion depends on the Newton polyhedra of $f_{1},\ldots
,f_{l} $ (see \cite{K}, \cite{O}). By \ constructing a log-principalization,
we give an explicit list \ for the possible poles of $Z_{\Phi }(s,%
\boldsymbol{f})$, $l\geq 1$, \ in the case in which $\boldsymbol{f}$ \ is
non degenerate with respect to $\Gamma \left( \boldsymbol{f}\right) $ (see
Theorem \ref{Theorem2a}). This theorem provides a generalization to the case
\ $l\geq 1$ of a well-known result that describes the poles of the local
zeta function associated to a non-degenerate polynomial in terms of the
corresponding Newton polyhedron \cite{D1a}, \cite{L-M}, \cite{Z1}. This
result was originally established by Varchenko \cite{Var}\ for local zeta
functions over $\mathbb{R}$. If $\boldsymbol{f}$ is non-degenerate with
respect to $\Gamma \left( \boldsymbol{f}\right) $, then $\lambda \left( 
\mathcal{I}_{\boldsymbol{f}}\right) $ can be computed from $\Gamma \left( 
\boldsymbol{f}\right) $ in the classical way (see Corollary \ref{coro2b}).

By using our notion of non-degeneracy and toroidal geometry we give an
explicit formula for $Z(s,\boldsymbol{f})$ and $Z_{0}(s,\boldsymbol{f})$, $%
l\geq 1$. This formula generalizes one given by Denef and Hoornaert in the
case $l=1$ \cite[Theorem 4.2]{D-H}, and one given by the second author for
the local zeta function of a monomial mapping \cite[Theorem 6.1]{Z1}.

The authors wish to thank the referee for his/her constructive remarks about
the paper. The first author would like to thank Robert Lazarsfeld for
suggestions and several inspiring conversations on poles of zeta functions,
and Orlando Villamayor for very useful information concerning
principalization in the analytic setting.

\section{The Igusa local zeta function of a polynomial mapping}

\subsection{Log-principalization and poles of the Igusa local zeta function}

We state the two versions of log-principalization of ideals that we will use
in this paper. The first is the `classical' algebraic formulation, see for
example \cite{E-N-V}, \cite{H}, \cite{W}. The second is in the context of $%
p- $adic analytic functions. It follows from \ the results in \cite{E-N-V},
see 5.11 in that paper (noticing that `Property D' there is valid in the $p-$%
adic analytic setting).

\begin{theorem}[Hironaka]
\label{ThHironaka}Let \ $X_{0}$ be a smooth algebraic variety over a field
of characteristic zero, and $\mathcal{I}$ a sheaf of ideals on $X_{0}$.
There exists a log-principalization \ of $\mathcal{I}$, that is a sequence 
\begin{equation*}
X_{0}\overset{\sigma _{1}}{\longleftarrow }X_{1}\overset{\sigma _{2}}{%
\longleftarrow }X_{2}\ldots \overset{\sigma _{i}}{\longleftarrow }X_{i}%
\overset{}{\longleftarrow }\ldots \overset{\sigma _{r}}{\longleftarrow }%
X_{r}=X
\end{equation*}%
of blow-ups $\sigma _{i}:X_{i-1}\longleftarrow X_{i}$ in smooth centers $%
C_{i-1}\subset X_{i-1}$ such that

\noindent (1) the exceptional divisor $E_{i}$ of the induced morphism $%
\sigma ^{i}=\sigma _{1}\circ \ldots \circ \sigma _{i}:X_{i}\longrightarrow
X_{0}$\ has only simple normal crossings and $C_{i}$ has simple normal
crossings with $E_{i}$, and

\noindent (2) the total transform $\left( \sigma ^{r}\right) ^{\ast }\left( 
\mathcal{I}\right) $ is the ideal of a simple normal crossings divisor $%
E^{\#}$. If the subscheme \ determined by $\mathcal{I}$ has no components \
of codimension one, then\ $E^{\#}$ is a natural combination of the
irreducible components of the divisor $E_{r}$.
\end{theorem}

\begin{remark}
We use notations like $(\sigma ^{r})^{\ast }(\mathcal{I})$ as in \cite{W}.
However, other authors use the notation $\mathcal{I}\mathcal{O}_{X}$ for the
same object, for example in \cite{E-N-V}. As many other authors we use the
term `log-principalization'. The terms `principalization'\ and
`monomialization' are also used for the same purpose by other authors.
\end{remark}

\begin{theorem}[\protect\cite{E-N-V}]
\label{thresolsing}Let $K$ be a $p-$adic field and $U$ an open submanifold
of $K^{n}$. Let $f_{1},\ldots ,f_{l}$ be $K-$analytic functions on $U$ such
that the ideal $\mathcal{I}_{\boldsymbol{f}}=\left( f_{1},\ldots
,f_{l}\right) $ is not trivial. Then there exists a log-principalization $%
\sigma :X_{K}\rightarrow U$ \ of $\mathcal{I}_{\boldsymbol{f}}$, that is,

\noindent (1) $X_{K}$ is an $n-$dimensional $K-$analytic manifold, $\sigma $
is a proper $K-$analytic map which is a composition of a finite number of
blow-ups in closed submanifolds, and which is an isomorphism outside of the
common zero set $Z_{K}$ of $f_{1},\ldots ,f_{l}$;

\noindent (2) $\sigma ^{-1}\left( Z_{K}\right) =\cup _{i\in T}E_{i}$, where
\ the $E_{i}$\ are closed submanifolds of $X_{K}$ of codimension one, each
equipped with a pair of positive integers $\left( N_{i},v_{i}\right) $
satisfying the following. At every \ point $b$ of $X_{K}$ there exist local
coordinates $\left( y_{1},\ldots ,y_{n}\right) $ on $X_{K}$ around $b$ such
that, if $E_{1},\ldots ,E_{p}$ are the $E_{i}$ containing $b$, we have \ on
some neighborhood of $b$ that $E_{i}$ is given by $y_{i}=0$ for $i=1,\ldots
,p$,

\begin{equation*}
\sigma ^{\ast }\left( \mathcal{I}_{\boldsymbol{f}}\right) \text{ is
generated by }\varepsilon \left( y\right)
\dprod\limits_{i=1}^{p}y_{i}^{N_{i}},
\end{equation*}%
and 
\begin{equation*}
\sigma ^{\ast }\left( dx_{1}\wedge \ldots \wedge dx_{n}\right) =\eta \left(
y\right) \left( \dprod\limits_{i=1}^{p}y_{i}^{v_{i}-1}\right) dy_{1}\wedge
\ldots \wedge dy_{n},
\end{equation*}%
where \ $\varepsilon \left( y\right) $, $\eta \left( y\right) $ are units in
the local ring of $X_{K}$\ at $b$.
\end{theorem}

The $\left( N_{i},v_{i}\right) $, $i\in T$, are called \textit{the numerical
data of} $\sigma $.

\bigskip

Let $K$ be a $p-$adic field. Let $f_{1},\ldots ,f_{l}$ be polynomials over $%
K $ or $K-$analytic functions on $U$ $\subset K^{n}$. We set $\mathcal{I}_{%
\boldsymbol{f}}$ to be the $K-$analytic ideal generated by \ the $f_{i}$; we
suppose it is not trivial. Let $\Phi :K^{n}\rightarrow \mathbb{C}$ \ or $%
U\rightarrow \mathbb{C}$ be \ a Schwartz-Bruhat function, that is, a locally
constant function with compact support. We associate to $\boldsymbol{f=}%
\left( f_{1},\ldots ,f_{l}\right) $ and $\Phi $\ the Igusa zeta function $%
Z_{\Phi }(s,\boldsymbol{f})$ as in the introduction. \ The following theorem
\ yields a new proof of its meromorphic continuation, but especially it
gives a list of its possible poles in terms of the numerical data of a
log-principalization.

\begin{theorem}
\label{Theorem1}The local zeta function $Z_{\Phi }(s,\boldsymbol{f})$ admits
a meromorphic continuation to the complex plane as a rational function of $%
q^{-s}$. Furthermore, the poles have the form 
\begin{equation*}
s=-\frac{v_{i}}{N_{i}}-\frac{2\pi \sqrt{-1}}{\log q}\frac{k}{N_{i}},\text{ \ 
}k\in \mathbb{Z}\text{,}
\end{equation*}%
where the $\left( N_{i},v_{i}\right) $ are the numerical data of a
log-principalization $\sigma :X_{K}\longrightarrow U$ of the ideal $\mathcal{%
I}_{\boldsymbol{f}}=\left( f_{1},\ldots ,f_{l}\right) $.
\end{theorem}

\begin{proof}
We pick a log-principalization $\sigma $ of $\mathcal{I}_{\boldsymbol{f}}$ \
as in Theorem \ref{thresolsing} and we use all notations that were
introduced there.

At every point $b\in X_{K}$ we can take \ a chart $\left( V,\phi _{V}\right) 
$ with coordinates $\left( y_{1},\ldots ,y_{n}\right) $, which may be
schrinked later when necessary. Let $g(y)$ be a generator of $\sigma ^{\ast
}\left( \mathcal{I}_{\boldsymbol{f}}\right) =\sigma ^{\ast }\left(
f_{1},\ldots ,f_{l}\right) $\ in $V$. Then%
\begin{equation*}
g(y)=\varepsilon (y)\dprod\limits_{i=1}^{p}y_{i}^{N_{i}},
\end{equation*}

\begin{equation*}
\sigma ^{\ast }\left( dx_{1}\wedge \ldots \wedge dx_{n}\right) =\eta
(y)\left( \dprod\limits_{i=1}^{p}y_{i}^{v_{i}-1}\right) dy_{1}\wedge \ldots
\wedge dy_{n},
\end{equation*}%
where $\varepsilon (y)$ and $\eta (y)$ are units of the local ring of $X_{K}$
at $b$. Furthermore, since $\sigma ^{\ast }\left( \mathcal{I}_{\boldsymbol{f}%
}\right) $\ is locally generated by $g(y)$ we have%
\begin{equation*}
f_{i}^{\ast }(y)=g(y)\widetilde{f}_{i}(y),
\end{equation*}%
for $i=1,\ldots ,l$, $y\in V$, where each $\widetilde{f}_{i}(y)$ is an
analytic function on $V$. And, since $g(y)\in \sigma ^{\ast }\left( \mathcal{%
I}_{\boldsymbol{f}}\right) $, we also have $\
g(y)=\dsum\nolimits_{i=1}^{l}a_{i}(y)f_{i}^{\ast }(y)$, \ with $a_{i}(y)$ an
analytic function on $V$ for each $i$; therefore 
\begin{equation*}
1=\dsum\limits_{i=1}^{l}a_{i}(y)\widetilde{f}_{i}(y)\text{, for }y\in V\text{%
.}
\end{equation*}%
Then there exists at least one index $i_{0}$ such that $\widetilde{f}%
_{i_{0}}(b)\neq 0$, hence we may assume that $\widetilde{f}_{i_{0}}(y)\neq 0$
on $V$ and that%
\begin{equation*}
\left\Vert \left( f_{1}^{\ast }(y),\ldots ,f_{l}^{\ast }(y)\right)
\right\Vert _{K}^{s}=\left\Vert \left( \left( \widetilde{f}_{i}(y)\right)
_{i\notin H},\left( \widetilde{f}_{i}(b)\right) _{i\in H}\right) \right\Vert
_{K}^{s}\left\vert g(y)\right\vert _{K}^{s},
\end{equation*}%
for $y\in V$. Here $H\subseteq \left\{ 1,\ldots ,n\right\} $ such that $%
\widetilde{f}_{i}(b)\neq 0\Leftrightarrow i\in H$. We may further suppose
that 
\begin{equation*}
\left\Vert \left( \left( \widetilde{f}_{i}(y)\right) _{i\notin H},\left( 
\widetilde{f}_{i}(b)\right) _{i\in H}\right) \right\Vert _{K}^{s}=\left\Vert
\left( \widetilde{f}_{i}(b)\right) _{i\in H}\right\Vert _{K}^{s}
\end{equation*}%
on $V$. Since $\sigma $\ is proper, $\sigma ^{-1}\left( \text{supp}\left(
\Phi \right) \right) $ is compact open in $X_{K}$, hence we can express it
as a finite \ disjoint union of compact open sets $B_{\alpha }$ such that
each $B_{\alpha }$ is contained in some $V$ above. Since $\Phi $\ is locally
constant we may assume (after \ subdividing $B_{\alpha }$) that $\left( \Phi
\circ \sigma \right) \mid _{B_{\alpha }}=$ $\left( \Phi \circ \sigma \right)
\left( b\right) $, $\left\vert \varepsilon \right\vert _{K}\mid _{B_{\alpha
}}=\left\vert \varepsilon \left( b\right) \right\vert _{K}$, $\left\vert
\eta \right\vert _{K}\mid _{B_{\alpha }}=\left\vert \eta \left( b\right)
\right\vert _{K}$, and $\phi _{V}\left( B_{\alpha }\right) =c+\pi
^{e_{0}}R_{K}^{n}$.

Denote by $D_{K}=\left( \text{div}\left( \sigma ^{\ast }\left( \mathcal{I}_{%
\boldsymbol{f}}\right) \right) \right) _{K}$. Since $\sigma :X_{K}\backslash
D_{K}\longrightarrow U\backslash \sigma \left( D_{K}\right) $ is bianalytic,
and \ $D_{K}$ has \ measure zero, we have 
\begin{equation*}
Z_{\Phi }(s,\boldsymbol{f})=\int\limits_{U\setminus \sigma \left(
D_{K}\right) }\Phi \left( x\right) \left\Vert \boldsymbol{f}(x)\right\Vert
_{K}^{s}\mid dx\mid
\end{equation*}%
\begin{equation*}
=\dsum\limits_{\alpha }\left( \Phi \circ \sigma \right) \left( b\right)
\left\vert \varepsilon \left( b\right) \right\vert _{K}^{s}\left\vert \eta
\left( b\right) \right\vert _{K}\left\Vert \left( \widetilde{f}%
_{i}(b)\right) _{i\in H}\right\Vert _{K}^{s}\dint\limits_{c+\pi
^{e_{0}}R_{K}^{n}}\dprod\limits_{1\leq i\leq p}\left\vert y_{i}\right\vert
^{N_{i}s+v_{i}-1}\left\vert dy\right\vert .
\end{equation*}%
The conclusion is now obtained by computing the integral in the previous
expression \ like in the case $l=1$ (see \cite[Lemma 8.2.1]{I2}).
\end{proof}

\begin{example}
\label{Ex1}Let $K$ be a $p-$adic field, and let $f_{1}(x,y)=y^{a}-x^{b}$, $%
f_{2}(x,y)=x^{a}-y^{b}$, with $a\,<b$, and for $j=3,\ldots ,M$, $M\geq 3$, $%
f_{j}\left( x,y\right) =x^{n_{j}}y^{m_{j}}h_{j}\left( x,y\right) $, with $%
n_{j}$, $m_{j}\geq a$, and $h_{j}\left( x,y\right) \in K\left[ x,y\right] $.
Set $\boldsymbol{f}=\left( f_{1},f_{2},f_{3},\ldots ,f_{M}\right) $, and \ $%
I_{\boldsymbol{f}}=\left( f_{1},f_{2},f_{3},\ldots ,f_{M}\right) $. Let $%
\Phi $ \ be a Schwartz-Bruhat \ function \ whose support \ is contained \ in
a sufficiently small neighborhood of the origin. A log-principalization of
the ideal $I_{\boldsymbol{f}}$ (over a neighborhood of the origin) is
obtained by blowing-up the origin of $K^{2}$. There is only one exceptional
curve $E=\mathbb{P}^{1}(K)$ whose numerical datum is $(a,2)$, and therefore
the possible poles of $Z_{\Phi }\left( s,\boldsymbol{f}\right) $\ have real
part $\frac{-2}{a}$. In \cite{Z2} an algorithm for computing a \ list of
candidates for the poles of $Z_{\Phi }\left( s,\boldsymbol{f}\right) $ in
terms of the numerical data of an embedded resolution of the divisor $\cup
_{j=1}^{M}f_{j}^{-1}\left( 0\right) $ was given. Since the $f_{j}\left(
x,y\right) $ are arbitrary polynomials for $3\leq j\leq M$, the mentioned
algorithm gives in general a very long list of possible poles.
\end{example}

\subsection{The largest real part of the poles of the Igusa zeta function}

Let $U$ be \ a compact open subset \ of $K^{n}$ and let $\boldsymbol{f}%
=(f_{1},\ldots ,f_{l}):U\longrightarrow K^{l}$ be an analytic\textit{\ }%
mapping. Recall that $Z_{U}\left( s,\boldsymbol{f}\right)
=\tint\nolimits_{U}\left\Vert \boldsymbol{f}(x)\right\Vert
_{K}^{s}\left\vert dx\right\vert $. The following lemma is known by the
experts, however we did not find \ a suitable reference for it; for the sake
of completeness we include its proof here.

\begin{lemma}
(1) $Z_{U}\left( s,\boldsymbol{f}\right) $ has no pole in $s$, i.e. $%
Z_{U}\left( s,\boldsymbol{f}\right) $ is a Laurent polynomial in $q^{-s}$ if
and only if there is no $x\in U$ such that \ $f_{1}\left( x\right) =\ldots
=f_{l}\left( x\right) =0$.

\noindent (2) If $0\in U$ and $\boldsymbol{f}(0)=0$, i.e. $f_{1}\left(
0\right) =\ldots =f_{l}\left( 0\right) =0$, then $Z_{U}\left( s,\boldsymbol{f%
}\right) $ has at least one pole in $s$.
\end{lemma}

\begin{proof}
(1) We first note \ that rationality of $Z_{U}\left( s,\boldsymbol{f}\right) 
$ implies the equivalence of the conditions \textquotedblleft $Z_{U}\left( s,%
\boldsymbol{f}\right) $ has no pole in $s$\textquotedblright\ and
\textquotedblleft $Z_{U}\left( s,\boldsymbol{f}\right) $ is a Laurent
polynomial in $q^{-s}$.\textquotedblright

\noindent ($\Leftarrow $) Since $\boldsymbol{f}:U\longrightarrow K^{l}$ is
continuous, also $\left\Vert \boldsymbol{f}\right\Vert _{K}:U\longrightarrow
q^{\mathbb{Z}}\cup \left\{ 0\right\} $ is continuous. If $0$ does not belong
to the image of $\left\Vert \boldsymbol{f}\right\Vert _{K}$, then \ there
are only finitely many values in the image because $U$ is compact. So $%
\tint\nolimits_{U}\left\Vert \boldsymbol{f}(x)\right\Vert _{K}^{s}\left\vert
dx\right\vert $ is a Laurent polynomial in $q^{-s}$.

\noindent ($\Rightarrow $) If $x_{0}\in U$ with \ $f_{1}\left( x_{0}\right)
=\ldots =f_{l}\left( x_{0}\right) =0$, by using the continuity of $%
\left\Vert \boldsymbol{f}\right\Vert _{K}$, there exist infinitely many $i$
such that there exists $x_{i}\in U$ with $\left\Vert \boldsymbol{f}\left(
x_{i}\right) \right\Vert _{K}=q^{-i}$. Since $U$ is open we have for all
those $i$ that the Haar measure of \ the set%
\begin{equation*}
\left\{ x\in U\mid \left\Vert \boldsymbol{f}\left( x\right) \right\Vert
_{K}=q^{-i}\right\}
\end{equation*}%
is positive. Therefore 
\begin{equation*}
Z_{U}\left( s,\boldsymbol{f}\right) =\tsum\limits_{j}vol\left( \left\{ x\in
U\mid \left\Vert \boldsymbol{f}\left( x\right) \right\Vert
_{K}=q^{-j}\right\} \right) \text{ }q^{-sj}
\end{equation*}%
is not a Laurent polynomial in $q^{-s}$.

(2) The second part follows directly from the first one.
\end{proof}

\begin{theorem}
\label{theorem1c}Let $\boldsymbol{f}=(f_{1},\ldots ,f_{l}):U\longrightarrow
K^{l}$ be an analytic\textit{\ mapping} defined on a compact open
neighborhood of the origin $U$ such that $\boldsymbol{f}(0)=0$. We take a
log-principalization $\sigma :X_{K}\rightarrow U$ as in Theorem \ref%
{thresolsing} with numerical data $\left( N_{i},v_{i}\right) $, $i\in T$.
Let $\lambda :=\lambda \left( \mathcal{I}_{\boldsymbol{f}}\right) =\min_{i}%
\frac{v_{i}}{N_{i}}$. Then $-\lambda \left( \mathcal{I}_{\boldsymbol{f}%
}\right) $ is the real part of a pole of $Z_{U}\left( s,\boldsymbol{f}%
\right) $. In particular, $\lambda \left( \mathcal{I}_{\boldsymbol{f}%
}\right) $ depends only on $\mathcal{I}_{\boldsymbol{f}}$.
\end{theorem}

\begin{proof}
The proof will be achieved by establishing that $q^{\lambda }$ is the
radius\ of convergence \texttt{R} of $Z_{U}\left( s,\boldsymbol{f}\right) $
considered as a function \ in $q^{-s}$. Certainly \texttt{R}$\geq q^{\lambda
}$, since (by Theorem \ref{Theorem1}) the candidate poles closest to the
origin have modulus $q^{\lambda }$. We shall show that \texttt{R}$\leq
q^{\lambda }$ by proving a lower bound for the \ coefficients of $%
Z_{U}\left( q^{-s},\boldsymbol{f}\right) $, considered as power series in $%
q^{-s}$:%
\begin{equation*}
Z_{U}\left( q^{-s},\boldsymbol{f}\right) =\tsum\limits_{j}vol\left( \left\{
x\in U\mid \left\Vert \boldsymbol{f}\left( x\right) \right\Vert
_{K}=q^{-j}\right\} \right) \text{ }q^{-sj}.
\end{equation*}

Take a \textit{generic} point $b$ on a component $E_{r}$ with $\frac{v_{r}}{%
N_{r}}=\lambda $, and a small enough chart $B$ ($\subset X_{K}$) around $b$
with coordinates $\left( y_{1},\ldots ,y_{n}\right) $ such that

\begin{equation*}
\sigma ^{\ast }\left( \mathcal{I}_{\boldsymbol{f}}\right) \text{ is
generated by }\varepsilon \left( y\right) y_{1}^{N_{r}},
\end{equation*}%
and 
\begin{equation*}
\sigma ^{\ast }\left( dx_{1}\wedge \ldots \wedge dx_{n}\right) =\eta \left(
y\right) y_{1}^{v_{r}-1}dy_{1}\wedge \ldots \wedge dy_{n},
\end{equation*}%
on $B$, where \ $\left\vert \varepsilon \right\vert _{K}$ and $\left\vert
\eta \right\vert _{K}$ are constant (and nonzero) on $B$. After an eventual $%
K-$analytic coordinate change, we may assume furthermore that $B=R_{K}^{n}$.

\textbf{Claim.} \textit{For }$j$\textit{\ big enough and divisible by }$%
N_{r} $\textit{\ we have }%
\begin{equation*}
vol\left( \left\{ x\in U\mid \left\Vert \boldsymbol{f}\left( x\right)
\right\Vert _{K}=q^{-j}\right\} \right) \geq Cq^{-j\lambda },
\end{equation*}%
\textit{where }$C$\textit{\ is a positive constant.}

By the above claim we have 
\begin{equation*}
\limsup_{i\rightarrow \infty }\left[ vol\left( \left\{ x\in U\mid \left\Vert 
\boldsymbol{f}\left( x\right) \right\Vert _{K}=q^{-i}\right\} \right) \right]
^{1/i}\geq q^{-\lambda }
\end{equation*}%
and hence 
\begin{equation*}
\mathtt{R}=\frac{1}{\limsup_{i\rightarrow \infty }\left[ vol\left( \left\{
x\in U\mid \left\Vert \boldsymbol{f}\left( x\right) \right\Vert
_{K}=q^{-i}\right\} \right) \right] ^{1/i}}\leq q^{\lambda }.
\end{equation*}

Therefore, since $Z_{U}\left( q^{-s},\boldsymbol{f}\right) $ is a \ rational
function of $q^{-s}$, we conclude that $uq^{\lambda }$ is a pole of $%
Z_{U}\left( q^{-s},\boldsymbol{f}\right) $, for some complex $N_{r}$-th root
of the unity $u$.

\textbf{Proof of the claim.} By the $p-$adic change of variables formula 
\cite[Proposition 7.4.1]{I2} we have ($B\subset \sigma ^{-1}\left( U\right) $%
):

\begin{equation*}
vol\left( \left\{ x\in U\mid \left\Vert \boldsymbol{f}\left( x\right)
\right\Vert _{K}=q^{-j}\right\} \right) \geq
\end{equation*}%
\begin{equation}
vol\left( \left\{ y\in B\mid \left\Vert \boldsymbol{f}\circ \sigma \left(
y\right) \right\Vert _{K}=q^{-j}\right\} \right) \cdot \left\vert \left( Jac%
\text{ }\sigma \right) \left( y\right) \right\vert _{K},  \label{cond1a}
\end{equation}%
where $Jac$ $\sigma $ is the Jacobian determinant of $\sigma $. With the
same reasoning as in the proof of Theorem \ref{Theorem1} we have that $%
\left\Vert \boldsymbol{f}\circ \sigma \left( y\right) \right\Vert
_{K}=C_{1}\left\vert \varepsilon \right\vert _{K}\left\vert y_{1}\right\vert
_{K}^{N_{r}}$ on $B$, where $C_{1}$ is a positive constant. So on $B$ we
have $\left\Vert \boldsymbol{f}\circ \sigma \left( y\right) \right\Vert
_{K}=q^{-j}$ if and only if $\left\vert y_{1}\right\vert
_{K}=C_{2}q^{-j/N_{r}}$, where $C_{2}$ is a positive constant. Hence%
\begin{equation}
vol\left( \left\{ y\in B\mid \left\Vert \boldsymbol{f}\circ \sigma \left(
y\right) \right\Vert _{K}=q^{-j}\right\} \right) =\left( 1-q^{-1}\right)
C_{2}q^{-j/N_{r}}.  \label{cond2}
\end{equation}%
Note that on this subset of $B$ we have%
\begin{equation}
\left\vert \left( Jac\text{ }\sigma \right) \left( y\right) \right\vert
_{K}=\left\vert \eta \right\vert _{K}\left\vert y_{1}\right\vert
_{K}^{v_{r}-1}=\left\vert \eta \right\vert _{K}C_{2}^{v_{r}-1}q^{-j\left(
v_{r}-1\right) /N_{r}}.  \label{cond3}
\end{equation}%
Combining (\ref{cond1a}), (\ref{cond2}) and (\ref{cond3}) yields 
\begin{equation*}
vol\left( \left\{ x\in U\mid \left\Vert \boldsymbol{f}\left( x\right)
\right\Vert _{K}=q^{-j}\right\} \right) \geq Cq^{-\lambda j},
\end{equation*}%
for some positive constant $C$.
\end{proof}

\begin{remark}
\label{rem0}(1) In \cite{I1} Igusa showed in the case $l=1$ that $-\lambda
\left( \mathcal{I}_{\boldsymbol{f}}\right) $ is a pole of $Z_{U}\left( s,%
\boldsymbol{f}\right) $ for a suitable compact open set $U$ containing the
origin. The argument uses Langlands' description of residues in terms of
principal value integrals \cite{Lan}. Furthermore, this argument is valid
for archimedean and non-archimedean local zeta functions (see also \cite[Th%
\'{e}or\`{e}me 5, part 3a, page 186]{AVG}, \cite{Var}).

\noindent (2) We note that $\lambda (\mathcal{I}_{\boldsymbol{f}})\geq lct(%
\mathcal{I}_{\boldsymbol{f}})$, where $lct(\mathcal{I}_{\boldsymbol{f}})$ is
the `log-canonical threshold' of $\mathcal{I}_{\boldsymbol{f}}$. This
well-known important invariant (see e.g. \cite{Ko}, \cite{MU}) is defined
analogously as $\lambda (\mathcal{I}_{\boldsymbol{f}})$ but in a geometric
setting, i.e. working over an algebraic closure of $K$. In order to obtain a
log-principalization in this context maybe more exceptional components are
needed, and then the inequality above could be strict.
\end{remark}

\subsubsection{Number of solutions of polynomial congruences}

Suppose that $f_{i}(x)$, $i=1,...,l$, are polynomials with coefficients in $%
R_{K}$. Let $N_{j}(\boldsymbol{f})$ be the number of solutions of $%
f_{i}(x)\equiv 0\,\ $mod$\,\ P_{K}^{j},\ i=1,...,l$, in $\left(
R_{K}/P_{K}^{j}\right) ^{n}$,$\,\ $and let \ $P(t,\boldsymbol{f})$ be the
series $\sum_{j=0}^{\infty }N_{j}(\boldsymbol{f})(q^{-n}t)^{j}$. The Poincar%
\'{e} series $P(t,\boldsymbol{f})$ is related to $Z(s,\boldsymbol{f})$ by
the formula $P(t,\boldsymbol{f})=\frac{1-tZ(s,\boldsymbol{f})}{1-t},$ \ $%
t=q^{-s},$ (cf. \cite[Theorem 2]{M2}). In the proof of the previous theorem
was established that $q^{\lambda }$ is the radius \ of convergence \texttt{R}
of $Z\left( s,\boldsymbol{f}\right) $ considered as a function \ in $q^{-s}$%
. By using this fact, \ and \ the above-mentioned relation between $P(t,%
\boldsymbol{f})$ and $Z(s,\boldsymbol{f})$, we obtain the following
corollary.

\begin{corollary}
\label{coro1}With the above notation, 
\begin{equation*}
\limsup_{j\rightarrow \infty }\left[ N_{j}(\boldsymbol{f})q^{-nj}\right] ^{%
\frac{1}{j}}=q^{-\lambda \left( \mathcal{I}_{\boldsymbol{f}}\right) },
\end{equation*}%
where\ $\lambda \left( \mathcal{I}_{\boldsymbol{f}}\right) =\min \left\{ 
\frac{v_{i}}{N_{i}}\right\} $, where $\left( N_{i},v_{i}\right) $ runs
through the numerical data of a log-principalization $\sigma
:X_{K}\longrightarrow R_{K}^{n}$ of the ideal $\mathcal{I}_{\boldsymbol{f}%
}=\left( f_{1},\ldots ,f_{l}\right) $.
\end{corollary}

Let $d$ be the maximal order of the poles of $P(t,\boldsymbol{f})$ with
modulus $q^{\lambda \left( \mathcal{I}_{\boldsymbol{f}}\right) }$. As a
consequence of the above corollary and of the rationality of $P(t,%
\boldsymbol{f})$ we have that $N_{j}(\boldsymbol{f})\leq
Cj^{d-1}q^{(n-\lambda \left( \mathcal{I}_{\boldsymbol{f}}\right) )j}$ for $j$
big enough, where $C$ is a positive constant. And by Remark \ref{rem0} (2),
we have then that $N_{j}(\boldsymbol{f})\leq Cj^{d-1}q^{(n-lct(\mathcal{I}_{%
\boldsymbol{f}}))j}$ for $j$ big enough.

\subsection{Denef's explicit formula}

For polynomials $f_{1},\ldots ,f_{l}$ over a number field $F$, we can
consider \ local zeta functions $Z_{W}(s,\boldsymbol{f,}K)$ for all
(non-archimedean) completions $K$ of $F$. When $l=1$, \ Denef presented in 
\cite[Theorem 3.1]{D2}\ an explicit formula, which is valid simultaneously
for almost all these zeta functions. \ His arguments extend to the several
polynomials case, by replacing resolution by log-principalization (as in
Theorem \ref{ThHironaka}).

\begin{theorem}
\label{Theorem1a}Let $F$ be a number field and $f_{i}(x)\in F\left[
x_{1},\ldots ,x_{n}\right] \ $\ for $i=1,\ldots ,l$. Let $\sigma
:X\rightarrow \mathbb{A}^{n}$ be a log-principalization of \ $I_{\boldsymbol{%
f}}=\left( f_{1},\ldots ,f_{l}\right) $ over $F$ as in Theorem \ref%
{ThHironaka}. Denote div$\left( \sigma ^{\ast }\left( \mathcal{I}_{%
\boldsymbol{f}}\right) \right) =\tsum\nolimits_{i\in T}N_{i}E_{i}$, and div$%
\left( \sigma ^{\ast }\left( dx_{1}\wedge \ldots \wedge dx_{n}\right)
\right) =\tsum\nolimits_{i\in T}\left( v_{i}-1\right) E_{i},$ where $E_{i}$, 
$i\in T$, are the irreducible components of the simple normal crossings
divisor given by the principal ideal $\sigma ^{\ast }\left( \mathcal{I}_{%
\boldsymbol{f}}\right) $. For every maximal ideal $P$ of the ring of
integers of $F$, we consider the completion $K$ of $F$ with respect to $P$.
Denote the valuation ring and the residue field of $K$ by $R$ and $\overline{%
K}$ $=\mathbb{F}_{q}$ respectively. Then for almost all completions $K$
(i.e. for all except a finite number) we have 
\begin{equation*}
Z_{W}\left( s,\boldsymbol{f},K\right) =q^{-n}\dsum\limits_{I\subseteq
T}c_{I}\dprod\limits_{i\in I}\frac{\left( q-1\right) q^{-N_{i}s-v_{i}}}{%
1-q^{-N_{i}s-v_{i}}},
\end{equation*}%
where $W\subset R^{n}$ is a union of cosets $mod$ $\left( P\right) ^{n}$,
and 
\begin{equation*}
c_{I}=\text{card}\left\{ a\in \overline{X}\left( \overline{K}\right) \mid
a\in \overline{E_{i}}\left( \overline{K}\right) \Leftrightarrow i\in I\text{%
; and }\overline{\sigma }(a)\in \overline{W}\right\} .
\end{equation*}%
Here $\overline{\cdot }$ denotes the reduction mod $P$, for which we refer
to \cite[Sect. 2]{D2}.
\end{theorem}

\begin{example}
Take $f_{1},f_{2},f_{3},\ldots ,f_{M}$\ as in Example \ref{Ex1} as being
defined over a number field $F$. Then the formula of Theorem \ref{Theorem1a}
\ for $W=\left( P\right) ^{2}$ yields 
\begin{equation*}
Z_{0}\left( s,\boldsymbol{f},K\right) =q^{-2}\left( q+1\right) \frac{\left(
q-1\right) q^{-as-2}}{1-q^{-as-2}}=\frac{\left( 1-q^{-2}\right) q^{-as-2}}{%
1-q^{-as-2}}.
\end{equation*}
\end{example}

\begin{example}
\label{Exam3}Let $K=\mathbb{Q}_{p}$, $f_{1}(x,y)=x$, $f_{2}(x,y)=x+p_{0}y$,
where $p_{0}$ is a fixed prime number, and let $\boldsymbol{f}=(f_{1},f_{2})$%
. A direct calculation shows that%
\begin{equation*}
Z(s,\boldsymbol{f},K)=\left\{ 
\begin{array}{ll}
\frac{1-p^{-2}}{1-p^{-2-s}}, & p\neq p_{0}, \\ 
&  \\ 
\frac{\left( 1-p^{-1}\right) \left( 1+p^{-1-s}\right) }{1-p^{-2-s}}, & 
p=p_{0}.%
\end{array}%
\right.
\end{equation*}%
A log-principalization for the ideal $\mathcal{I}_{\boldsymbol{f}}$ is
attained \ by blowing-up the origin. One easily verifies that the expression
\ for $p\neq p_{0}$ is the one given by Theorem \ref{Theorem1a}.
\end{example}

As a consequence of Theorem \ref{Theorem1} (or \cite{D1}, \cite{M2}\ ) $%
Z_{W}(s,\boldsymbol{f})$ can be written as 
\begin{equation*}
Z_{W}(s,\boldsymbol{f})=\frac{P(T)}{Q(T)},
\end{equation*}%
where $P(T)$ and $Q(T)$ are polynomials in $T=q^{-s}$ with rational
coefficients. We define $\deg $ $Z_{W}(s,\boldsymbol{f})=\deg $ $P(T)$ $-$ $%
\deg $ $Q(T)$, where $\deg $ means `degree'.

\begin{corollary}
Let $f_{i}(x)\in F\left[ x_{1},\ldots ,x_{n}\right] $ for $i=1,\ldots ,l$.
For almost all completions $K$ of $F$ we have \ $\deg $ $Z(s,\boldsymbol{f}%
,K)\leq 0$ and $\deg $ $Z_{0}(s,\boldsymbol{f},K)=0$. Moreover if all $f_{i}$
are homogeneous of degree $d$, then $\deg Z(s,\boldsymbol{f},K)=-d$.
\end{corollary}

The proof \ follows from the explicit formula (Theorem \ref{Theorem1a}) by
analogous arguments \ as in \cite{D2} (or \ \cite{I2}) where the case $l=1$
is treated. We should mention that by using model-theoretic arguments Denef
already showed \ the above result (see \cite[Theorem 5.2, and Example 5.4]%
{D2}). So in this paper we give a geometric proof of this fact.

Note that for the case $p=p_{0}$ in Example \ref{Exam3} it is not true that $%
\deg Z(s,\boldsymbol{f},\mathbb{Q}_{p})=-1$, though $f_{1}$, $f_{2}$ are
homogeneous of degree $1$.

\begin{example}
\label{Exam4}Let $\boldsymbol{f}=\left( f_{1},f_{2}\right) =\left(
x^{3}-xy,y\right) $. One easily constructs a log-principalization of the
ideal $\mathcal{I}_{\boldsymbol{f}}=\left( x^{3}-xy,y\right) $ as a
composition of three blow-ups. The numerical data of the three \ exceptional
components in $\sigma ^{-1}\left( \text{supp }\mathcal{I}_{\boldsymbol{f}%
}\right) =\sigma ^{-1}\left( 0\right) $ are $\left( 1,2\right) $, $\left(
2,3\right) $, $\left( 3,4\right) $ respectively. \ So Theorem \ref{Theorem1}
yields $-2$, $-3/2$, $-4/3$ as possible (real parts of) candidate poles of $%
Z(s,\boldsymbol{f})$. However, in the formula \ of Theorem \ref{Theorem1a}
the first two candidate poles cancel:

\begin{equation*}
Z(s,\boldsymbol{f})=q^{-2}\{\left( q^{2}-1\right) +q\frac{\left( q-1\right)
q^{-2-s}}{1-q^{-2-s}}+\left( q-1\right) \frac{\left( q-1\right) q^{-3-2s}}{%
1-q^{-3-2s}}+q\frac{\left( q-1\right) q^{-4-3s}}{1-q^{-4-3s}}
\end{equation*}%
\begin{equation*}
+\frac{\left( q-1\right) ^{2}q^{-5-3s}}{\left( 1-q^{-2-s}\right) \left(
1-q^{-3-2s}\right) }+\frac{\left( q-1\right) ^{2}q^{-7-5s}}{\left(
1-q^{-3-2s}\right) \left( 1-q^{-4-3s}\right) }\}
\end{equation*}%
\begin{equation*}
=q^{-2}\frac{q-1}{1-q^{-4-3s}}\left( q+1+q^{-1-s}+q^{-2-2s}\right) .
\end{equation*}%
We shall present an alternative formula to compute this example in Section
4, where only one candidate pole will appear.
\end{example}

\begin{example}
\label{Example5}Let $\boldsymbol{f}=\left( f_{1},f_{2}\right) =\left(
y^{2}-x^{3},y^{2}-z^{2}\right) $. We shall compute $Z_{0}\left( s,%
\boldsymbol{f}\right) $\ by means of a log-principalization of $\mathcal{I}_{%
\boldsymbol{f}}=\left( y^{2}-x^{3},y^{2}-z^{2}\right) $. Note that the
support of $\mathcal{I}_{\boldsymbol{f}}$ has two $1-$dimensional components 
$C$ and $C^{\prime }$ with a singularity at the origin of $K^{3}$.

We first blow up the origin yielding the exceptional surface $E_{1}\left(
\cong \mathbb{P}^{2}\right) $ with $\left( N_{1},v_{1}\right) =\left(
2,3\right) $. The strict transform of $C$ and $C^{\prime }$ and $E_{1}$ have
one common point. Next we blow up this point obtaining the new exceptional
surface $E_{2}\left( \cong \mathbb{P}^{2}\right) $ with $\left(
N_{2},v_{2}\right) =\left( 3,5\right) $. At this stage (the strict
transforms of) $C$ and $C^{\prime }$ are disjoint and both meet $E_{2}$ in
one point \ of the intersection of $E_{2}$ with (the strict transform of) $%
E_{1}$. Now we blow up the curve $E_{1}\cap E_{2}$; the new exceptional
component $E_{3}$ is a ruled surface over that curve and $%
(N_{3},v_{3})=(6,8) $. We have that $E_{3}\cap E_{1}$ and $E_{3}\cap E_{2}$
are disjoint sections of $E_{3}$, and $C$ and $C^{\prime }$ intersect $E_{3}$
transversely outside $E_{3}\cap E_{1}$ and $E_{3}\cap E_{2}$. Finally we
blow up $C$ and $C^{\prime }$, yielding the last two exceptional surfaces $%
E_{4}$ and $E_{4}^{\prime }$ with numerical data $(1,2)$. The formula of
Theorem \ref{Theorem1a} yields%
\begin{equation*}
Z_{0}(s,\mathbf{f})=q^{-3}\left( (q^{2}+q)\frac{(q-1)q^{-3-2s}}{1-q^{-3-2s}}%
+q^{2}\frac{(q-1)q^{-5-3s}}{1-q^{-5-3s}}\right.
\end{equation*}%
\begin{equation*}
+(q^{2}-3)\frac{(q-1)q^{-8-6s}}{1-q^{-8-6s}}+(q+1)\frac{(q-1)^{2}q^{-11-8s}}{%
(1-q^{-3-2s})(1-q^{-8-6s})}
\end{equation*}%
\begin{equation*}
\left. +(q+1)\frac{(q-1)^{2}q^{-13-9s}}{(1-q^{-5-3s})(1-q^{-8-6s})}+2(q+1)%
\frac{(q-1)^{2}q^{-10-7s}}{(1-q^{-2-s})(1-q^{-8-6s})}\right)
\end{equation*}%
\begin{equation*}
=q^{-3}(q-1)\frac{N\left( q^{-s}\right) }{(1-q^{-2-s})(1-q^{-8-6s})},
\end{equation*}%
where 
\begin{eqnarray*}
N\left( q^{-s}\right)
&=&(q^{2}-q-1)q^{-10-7s}+(q^{2}+q-1)q^{-8-6s}-(q+1)q^{-7-5s} \\
&&+q^{-4-4s}-q^{-4-3s}+(q+1)q^{-2-2s}.
\end{eqnarray*}

Note that the candidate poles $-3/2$ and $-5/3$ cancel.
\end{example}

\subsection{\label{motivic}Motivic and topological zeta functions}

The analogue of the original explicit formula of Denef plays an important
role in the study of the motivic \ zeta function associated to one regular
function \cite{D-L}. One can associate more generally a motivic zeta
function to any sheaf of ideals on a smooth variety, and obtain a similar
formula for it in terms of a log-principalization using the argument of \cite%
{D-L}. We just formulate \ the more general definition and formula,
referring to e.g. \cite{D-L2}, \cite{V1}\ for the notion of jets and
Grothendieck ring.

\begin{definition}
\label{defmotzeta}Let $Y$ be a smooth algebraic variety of dimension $n$
over over a field $F$ of characteristic zero, and $\mathcal{I}$ a sheaf of
ideals on $Y$. Let $W$ be a subvariety of $Y$. Denote for $i\in \mathbb{N}$
by $\mathfrak{X}_{i,W\text{ }}$ the variety of $i-$jets $\gamma $ on $Y$ \
with origin in $W$ for which $ord_{t}\left( \gamma ^{\ast }\mathcal{I}%
\right) =i$. The \textit{motivic zeta function} associated \ to $\mathcal{I}$
(and $W$) is the formal power series%
\begin{equation*}
Z_{W}\left( \mathcal{I},T\right) =\tsum\limits_{i\geq 0}\left[ \mathfrak{X}%
_{i,W\text{ }}\right] \left( \mathbb{L}^{-n}T\right) ^{i},
\end{equation*}%
where $\left[ \cdot \right] $ \ denotes the class of a variety in the
Grothendieck ring of algebraic varieties over $F$, and $\mathbb{L}=\left[ 
\mathbb{A}^{1}\right] $.
\end{definition}

\begin{theorem}
\label{Theorem2b}Let $\sigma :X\rightarrow Y$ be a log-principalization of $%
\mathcal{I}$. With the analogous notation $E_{i}$, $N_{i}$, $v_{i}$, $(i\in
T)$ \ as before, and also \ $E_{I}^{\circ }:=\left( \cap _{i\in
I}E_{i}\right) \setminus \left( \cup _{k\notin I}E_{k}\right) $ for $%
I\subset T$, we have%
\begin{equation*}
Z_{W}\left( \mathcal{I},T\right) =\tsum\limits_{I\subset T}\left[
E_{I}^{\circ }\cap \sigma ^{-1}W\right] \tprod\limits_{i\in I}\frac{\left( 
\mathbb{L}-1\right) T^{N_{i}}}{\mathbb{L}^{v_{i}}-T^{N_{i}}}.
\end{equation*}%
In particular $Z_{W}\left( \mathcal{I},T\right) $ is rational in $T$.
\end{theorem}

Specializing to topological Euler characteristics, denoted by $\chi \left(
\cdot \right) $, as in \cite[(2.3)]{D-L} or \cite[(6.6)]{V1} we obtain the
expression 
\begin{equation*}
Z_{top,W}\left( \mathcal{I},s\right) :=\tsum\limits_{I\subset T}\chi \left(
E_{I}^{\circ }\cap \sigma ^{-1}W\right) \tprod\limits_{i\in I}\frac{1}{%
v_{i}+N_{i}s}\in \mathbb{Q}(s),
\end{equation*}%
which is then independent of the chosen log-principalization. (When the base
field is not the complex numbers, we consider $\chi (\cdot )$ in \'{e}tale $%
\overline{\mathbb{Q}}$- cohomology as in [8].) It can be taken as a
definition for the \textit{topological zeta function} associated to $%
\mathcal{I}$ (and $W$), generalizing the original one of Denef and Loeser
associated to one polynomial \cite{D-L3}.

\section{Newton polyhedra and non-degeneracy conditions}

\subsection{Newton polyhedra}

We set $\mathbb{R}_{+}:=\{x\in \mathbb{R}\mid x\geqslant 0\}$.

Let $G$ be a nonempty subset of $\mathbb{N}^{n}$. The \textit{Newton
polyhedron }$\Gamma =\Gamma \left( G\right) $ associated to $G$ is the
convex hull in $\mathbb{R}_{+}^{n}$ \ \ of the set $\cup _{m\in G}\left( m+%
\mathbb{R}_{+}^{n}\right) $. For instance classically one associates a 
\textit{Newton polyhedron (at the origin) to } $g(x)=\sum_{m}c_{m}x^{m}$ ($%
x=\left( x_{1},\ldots ,x_{n}\right) $, $g(0)=0$), \ being a nonconstant
polynomial function over $K$ or $K-$analytic function \ in a neighborhood of
the origin, where $G=$supp$(g)$ $:=$ $\left\{ m\in \mathbb{N}^{n}\mid
c_{m}\neq 0\right\} $.\ Further we will associate more generally a Newton
polyhedron to an analytic mapping.

We fix a Newton polyhedron $\Gamma $\ as above. We first collect some
notions and results about Newton polyhedra that \ will be used in the next
sections. Let $\left\langle \cdot ,\cdot \right\rangle $ denote the usual
inner product of $\mathbb{R}^{n}$, and identify \ the dual space of $\mathbb{%
R}^{n}$ with $\mathbb{R}^{n}$ itself by means of it.

For $a\in \mathbb{R}_{+}^{n}$, we define 
\begin{equation*}
d(a,\Gamma )=d(a)=\min_{x\in \Gamma }\left\langle a,x\right\rangle ,
\end{equation*}%
and \textit{the first meet locus }$F(a)$ of $a$ as \ 
\begin{equation*}
F(a):=\{x\in \Gamma \mid \left\langle a,x\right\rangle =d(a)\}.
\end{equation*}%
The first meet locus \ is a face of $\Gamma $. Moreover, if $a\neq 0$, $F(a)$
is a proper face of $\Gamma $.

We define an equivalence relation in \ $\mathbb{R}_{+}^{n}$ by taking $a\sim
a^{\prime }\Leftrightarrow F(a)=F(a^{\prime })$. The equivalence classes of
\ $\sim $ are sets of the form 
\begin{equation*}
\Delta _{\tau }=\{a\in \mathbb{R}_{+}^{n}\mid F(a)=\tau \},
\end{equation*}%
where $\tau $ \ is a face of $\Gamma $.

We recall that the cone strictly\ spanned \ by the vectors $a_{1},\ldots
,a_{r}\in \mathbb{R}_{+}^{n}\setminus \left\{ 0\right\} $ is the set $\Delta
=\left\{ \lambda _{1}a_{1}+...+\lambda _{r}a_{r}\mid \lambda _{i}\in \mathbb{%
R}_{+}\text{, }\lambda _{i}>0\right\} $. If $a_{1},\ldots ,a_{r}$ are
linearly independent over $\mathbb{R}$, $\Delta $ \ is called \ a \textit{%
simplicial cone}. \ If \ $a_{1},\ldots ,a_{r}\in \mathbb{Z}^{n}$, we say $%
\Delta $\ is \ a \textit{rational cone}. If $\left\{ a_{1},\ldots
,a_{r}\right\} $ is a subset of a basis \ of the $\mathbb{Z}$-module $%
\mathbb{Z}^{n}$, we call $\Delta $ a \textit{simple cone}.

A precise description of the geometry of the equivalence classes modulo $%
\sim $ is as follows. Each \textit{facet} (i.e. a face of codimension one)\ $%
\gamma $ of $\Gamma $\ has a unique vector $a(\gamma )=(a_{\gamma ,1},\ldots
,a_{\gamma ,n})\in \mathbb{N}^{n}\mathbb{\setminus }\left\{ 0\right\} $, \
whose nonzero coordinates are relatively prime, which is perpendicular to $%
\gamma $. We denote by $\mathfrak{D}(\Gamma )$ the set of such vectors. The
equivalence classes are rational cones of the form 
\begin{equation*}
\Delta _{\tau }=\{\sum\limits_{i=1}^{r}\lambda _{i}a(\gamma _{i})\mid
\lambda _{i}\in \mathbb{R}_{+}\text{, }\lambda _{i}>0\},
\end{equation*}%
where $\tau $ runs through the set of faces of $\Gamma $, and $\gamma _{i}$, 
$i=1,\ldots ,r$\ are the facets containing $\tau $. We note that $\Delta
_{\tau }=\{0\}$ if and only if $\tau =\Gamma $. The family $\left\{ \Delta
_{\tau }\right\} _{\tau }$, with $\tau $ running over\ the proper faces of $%
\Gamma $, is a partition of $\mathbb{R}_{+}^{n}\backslash \{0\}$; we call
this partition a \textit{\ polyhedral subdivision of }\ $\mathbb{R}_{+}^{n}$%
\ \textit{subordinated} to $\Gamma $. We call $\left\{ \overline{\Delta }%
_{\tau }\right\} _{\tau }$, the family formed by the topological closures of
the $\Delta _{\tau }$, a \textit{\ fan} \textit{subordinated} to $\Gamma $.

Each cone $\Delta _{\tau }$\ can be partitioned \ into a finite number of
simplicial cones $\Delta _{\tau ,i}$. In addition, the subdivision can be
chosen such that each $\Delta _{\tau ,i}$ is spanned by part of $\mathfrak{D}%
(\Gamma )$ . Thus from the above considerations we have the following
partition of $\mathbb{R}_{+}^{n}\backslash \{0\}$:

\begin{equation}
\mathbb{R}_{+}^{n}\backslash \{0\}=\bigcup\limits_{\tau \text{ }}\left(
\bigcup\limits_{i=1}^{l_{\tau }}\Delta _{\tau ,i}\right) ,
\label{simplisubv}
\end{equation}%
where $\tau $ runs \ over the proper faces of $\Gamma $, and each $\Delta
_{\tau ,i}$ \ is a simplicial cone contained in $\Delta _{\tau }$.\ We will
say that $\left\{ \Delta _{\tau ,i}\right\} $ is a \textit{simplicial
polyhedral subdivision of }\ $\mathbb{R}_{+}^{n}$\ \textit{subordinated} to $%
\Gamma $; and that $\left\{ \overline{\Delta }_{\tau ,i}\right\} $ is a 
\textit{simplicial fan} \textit{subordinated} to $\Gamma $.

By adding new rays , each simplicial cone can be partitioned further into a
finite number of simple cones. In this way we obtain a \textit{simple
polyhedral subdivision} of $\mathbb{R}_{+}^{n}$\ \textit{subordinated} to $%
\Gamma $; and a \textit{simple fan} \textit{subordinated} to $\Gamma $ (see
e.g. \cite{K-M-S}).

\subsection{The Newton polyhedron associated to an analytic mapping}

Let $\boldsymbol{f}=(f_{1},\ldots ,f_{l})$, $\boldsymbol{f}\left( 0\right)
=0 $, be a nonconstant polynomial mapping, or more generally, an analytic
mapping defined \ on a neighborhood $U\subseteq K^{n}$ of the origin. In
this paper we associate to $\boldsymbol{f}$ a Newton polyhedron $\Gamma
\left( \boldsymbol{f}\right) :=\Gamma \left( \cup _{i=1}^{l}\text{supp}%
\left( f_{i}\right) \right) $, \ and a non-degeneracy condition to $%
\boldsymbol{f}$ and $\Gamma \left( \boldsymbol{f}\right) $.

If $f_{i}\left( x\right) =\tsum\nolimits_{m}c_{m,i}x^{m}$, and $\tau $ is a
face of $\Gamma \left( \boldsymbol{f}\right) $, we set%
\begin{equation*}
f_{i,\tau }\left( x\right) :=\tsum\limits_{m\in \text{supp}(f_{i})\cap \tau
}c_{m,i}x^{m}.
\end{equation*}

\begin{definition}
\label{def1}(1) Let $\boldsymbol{f}=(f_{1},\ldots ,f_{l}):U\longrightarrow
K^{l}$ be a nonconstant analytic mapping satisfying $\boldsymbol{f}\left(
0\right) =0$. The mapping \ $\boldsymbol{f}$ is called \textit{strongly
non-degenerate at the origin with respect to }$\Gamma (\boldsymbol{f})$,%
\textit{\ }if for \ any \textit{compact} face $\tau \subset \Gamma (%
\boldsymbol{f})$ and any $z\in \left\{ z\in \left( K^{\times }\right)
^{n}\mid f_{1,\tau }(z)=\ldots =f_{l,\tau }(z)=0\right\} $ it verifies that $%
rank_{K}\left[ \frac{\partial f_{i,_{\tau }}}{\partial x_{j}}\left( z\right) %
\right] =\min \{l,n\}$.

\noindent (2) Let $\boldsymbol{f}=(f_{1},\ldots ,f_{l}):K^{n}\longrightarrow
K^{l}$ be a nonconstant polynomial mapping satisfying $\boldsymbol{f}\left(
0\right) =0$. The mapping $\boldsymbol{f}$\ is called \textit{strongly} \ 
\textit{non-degenerate with \ respect to }$\Gamma (\boldsymbol{f})$,\textit{%
\ }if for any face $\ \tau \subset \Gamma (\boldsymbol{f})$, including $%
\Gamma (\boldsymbol{f})$ itself, and \ \ any $z\in \left\{ z\in \left(
K^{\times }\right) ^{n}\mid f_{1,\tau }(z)=\ldots =f_{l,\tau }(z)=0\right\} $
\ \ \ it verifies that $rank_{K}\left[ \frac{\partial f_{i,_{\tau }}}{%
\partial x_{j}}\left( z\right) \right] $ $=$ \ $\min \{l,n\}$.
\end{definition}

\begin{remark}
\label{rem1}Let $\boldsymbol{f}=(f_{1},\ldots ,f_{l}):U\longrightarrow K^{l}$
be a nonconstant analytic mapping satisfying $\boldsymbol{f}\left( 0\right)
=0$.

\noindent (1) Let $\gamma $ be a face of $\Gamma (\boldsymbol{f})$ for which
the rank condition in Definition \ref{def1} is satisfied. If supp$%
(f_{i})\cap \gamma \neq \varnothing \Leftrightarrow i\in I_{\gamma }$ for a
non-empty subset $I_{\gamma }$ $\subseteq \left\{ 1,\ldots ,l\right\} $
satisfying card$\left( I_{\gamma }\right) <\min \{l,n\}$, then necessarily%
\begin{equation*}
\dbigcap\limits_{i\in I_{\gamma }}\left\{ z\in \left( K^{\times }\right)
^{n}\mid f_{i,\gamma }\left( z\right) =0\right\} =\varnothing .
\end{equation*}

\noindent (2) If for a given face $\gamma $ at least one \ $f_{i,\gamma }$
is a monomial, then the rank condition \ on $\gamma $ is satisfied. This is
in particular true if $\gamma $ is a point.
\end{remark}

\begin{example}
\label{example1} Let $\boldsymbol{f}(x,y)=\left( x^{3}-xy,y\right) $. The
mapping $\boldsymbol{f}$ is strongly non-degenerate \ at the origin with
respect to $\Gamma (\boldsymbol{f})$, and also strongly non-degenerate with
respect to $\Gamma (\boldsymbol{f})$.
\end{example}

\begin{example}
\label{example3}Let $\boldsymbol{f}\left( x,y,z\right) =\left(
x^{2},y^{2},z^{2},xy,xz,yz\right) $. Then $\boldsymbol{f}$ is strongly
non-degenerate at the origin with respect to $\Gamma \left( \boldsymbol{f}%
\right) $, and also strongly non-degenerate with respect to $\Gamma \left( 
\boldsymbol{f}\right) $,
\end{example}

\subsubsection{Monomial mappings}

Any monomial mapping is strongly non-degenerate at the origin with respect
to its Newton polyhedron. If $\boldsymbol{f}_{0}$ is a fixed monomial
mapping with Newton polyhedron $\Gamma \left( \boldsymbol{f}_{0}\right) $,
and $\boldsymbol{f}=\boldsymbol{f}_{0}+\boldsymbol{g}$ is a deformation of \ 
$\boldsymbol{f}_{0}$ such that all the monomials in $\boldsymbol{g}$ have
exponents in the interior of $\Gamma \left( \boldsymbol{f}_{0}\right) $,
then $\boldsymbol{f}$ is \ strongly non-degenerate at the origin with \
respect to \ $\Gamma \left( \boldsymbol{f}\right) =\Gamma \left( \boldsymbol{%
f}_{0}\right) $. This type of mapping was introduced by the second author in 
\cite[Definition 6.1]{Z2}. Furthermore, the corresponding local zeta
function \ can be computed by using a simple polyhedral subdivision
subordinated to $\Gamma \left( \boldsymbol{f}_{0}\right) $ \cite[Theorem 6.1]%
{Z2}.

\subsubsection{Saia's non-degeneracy condition}

In \cite{S} Saia introduced \ the following notion of non-degeneracy for
ideals. Let $I=(f_{1},\ldots f_{l})$ \ be a polynomial ideal. $I$ is
non-degenerate with respect to $\Gamma \left( I\right) $ (where $\Gamma
\left( I\right) $ $=\Gamma \left( \cup _{i=1}^{l}\text{supp}\left(
f_{i}\right) \right) $), if for every compact face $\tau $ of $\Gamma \left(
I\right) $, the system of equations $f_{1,\tau }(z)=0,\ldots f_{l,\tau
}(z)=0 $ does not have a solution in the torus $\left( K^{\times }\right)
^{n}$. Thus \ Saia's notion of non-degeneracy is a particular case of our
notion of non-degeneracy. Saia's notion of non-degeneracy plays an important
role in \ the study of the integral closure of ideals.

\subsubsection{Khovanskii's non-degeneracy condition}

Now we discuss the relation \ between our notion of non-degeneracy and
Khovanskii's notion of non-degeneracy of an analytic mapping with respect to
several Newton polyhedra (\cite{K}, see also \cite{O}).\ Given a positive
vector $a$ (i.e. $a\in \left( \mathbb{N\setminus }\left\{ 0\right\} \right)
^{n}$), and an analytic mapping $g$, we set $g_{a}(x):=$ $g_{F(a)}(x)$,
where $F(a)$ is the first meet locus \ of $a$ with respect to $\Gamma \left(
g\right) $. To make explicit the dependence between $F(a)$ and $\Gamma
\left( g\right) $ we shall write $F(a,\Gamma \left( g\right) )$ instead of $%
F(a)$.

\begin{definition}
\label{Khovan}A nonconstant analytic mapping $\boldsymbol{f}=(f_{1},\ldots
,f_{l}):U\longrightarrow K^{l}$, $\boldsymbol{f}\left( 0\right) =0$, is
non-degenerate with respect to $(\Gamma \left( f_{1}\right) ,\ldots ,\Gamma
\left( f_{l}\right) )$, if for any \ positive vector $a$ and any $z\in
\left\{ z\in \left( K^{\times }\right) ^{n}\mid f_{1,a}(z)=\ldots
=f_{l,a}(z)=0\right\} $ it verifies that 
\begin{equation*}
rank_{K}\left[ \frac{\partial f_{i,a}}{\partial x_{j}}\left( z\right) \right]
=\min \{l,n\}.
\end{equation*}%
Here $f_{j,a}(z)=f_{j,F(a,\Gamma \left( f_{j}\right) )}(z)$ for every $j$.
\end{definition}

The above definition is equivalent \ to the non-degeneracy notion given by
Oka in \cite{O}, that is in turn a reformulation of the notion of
non-degeneracy introduced by Khovanskii in \cite{K}.

\begin{remark}
Let $\boldsymbol{f}=(f_{1},\ldots ,f_{l}):U\longrightarrow K^{l}$ be a
nonconstant analytic mapping satisfying $\boldsymbol{f}\left( 0\right) =0$.
Then $\Gamma \left( \boldsymbol{f}\right) $ is the convex hull in $\left( 
\mathbb{R}_{+}\right) ^{n}$ of $\ \cup _{j=1}^{l}\Gamma \left( f_{j}\right) $%
. \ This assertion follows from the fact \ that \ for any subsets $A$, $B$ $%
\subseteq $ $\left( \mathbb{R}_{+}\right) ^{n}$, $\overline{A\cup B}=%
\overline{\overline{A}\cup \overline{B}}$, where the bar denotes the convex
hull in $\left( \mathbb{R}_{+}\right) ^{n}$.
\end{remark}

The following is the relation between Khovanskii's non-degeneracy notion and
the one introduced here.

\begin{proposition}
Let $\boldsymbol{f}=(f_{1},\ldots ,f_{l}):U\longrightarrow K^{l}$ be an
analytic mapping strongly non-degenerate at the origin with respect to $%
\Gamma (\boldsymbol{f})$. Then $\boldsymbol{f}$\ is non-degenerate with
respect to\ 
\begin{equation*}
(\Gamma \left( f_{1}\right) ,\ldots ,\Gamma \left( f_{l}\right) ).
\end{equation*}
\end{proposition}

\begin{proof}
Let $a\in \left( \mathbb{N\setminus }\left\{ 0\right\} \right) ^{n}$ be a
fixed positive vector. We set $\Gamma =\Gamma (\boldsymbol{f})$, $\Gamma
_{j}=\Gamma \left( f_{j}\right) $, $j=1,\ldots ,l$. \ Since $\Gamma
_{j}\subseteq \Gamma $ by the above remark, 
\begin{equation*}
d(a,\Gamma )=\min_{x\in \Gamma }\left\langle a,x\right\rangle \leq
d(a,\Gamma _{j})=\min_{x\in \Gamma _{j}}\left\langle a,x\right\rangle ,
\end{equation*}%
for $j=1,\ldots ,l$. \ We define $I\subseteq \left\{ 1,\ldots ,l\right\} $ \
by the condition 
\begin{equation*}
j\in I\Leftrightarrow d(a,\Gamma )=d(a,\Gamma _{j}).
\end{equation*}%
Note that $I\neq \varnothing $. Then, \ if $\tau :=F(a,\Gamma )$,

\begin{equation*}
F(a,\Gamma _{j})\subseteq \tau \text{, for \ }j\in I,
\end{equation*}%
and%
\begin{equation}
f_{j,\tau }(x)=\left\{ 
\begin{array}{ll}
f_{j,a}(x), & \text{ }j\in I, \\ 
&  \\ 
0, & \text{ }j\in I^{c}.%
\end{array}%
\right.  \label{con1}
\end{equation}%
If card$(I)<\min \left\{ l,n\right\} $, then by \ Remark \ref{rem1} the
system of equations 
\begin{equation*}
f_{j,\tau }(x)=0\text{, }j\in I\text{, has no solutions in }\left( K^{\times
}\right) ^{n}.
\end{equation*}%
Hence by \ using (\ref{con1}) the system of equations%
\begin{equation*}
f_{j,a}(x)=0,\text{ }j=1,\ldots ,l,\text{ has no solutions in }\left(
K^{\times }\right) ^{n},
\end{equation*}%
and \ so the condition on $a$ in Definition \ref{Khovan} is satisfied.

Now, we may assume that card$(I)\geq \min \left\{ l,n\right\} $, and that 
\begin{equation*}
f_{j,\tau }(x)=0,\text{ }j\in I,\text{ has solutions in }\left( K^{\times
}\right) ^{n}.
\end{equation*}%
Since $\boldsymbol{f}$ is strongly non-degenerate with respect to $\Gamma
\left( \boldsymbol{f}\right) $, it follows that 
\begin{equation*}
rank_{K}\left[ \frac{\partial f_{j,\tau }}{\partial x_{i}}\left( z\right) %
\right] =rank_{K}\left[ \frac{\partial f_{j,\tau }}{\partial x_{i}}\left(
z\right) \right] _{\substack{ j\in I  \\ 1\leq i\leq n}}=\min \{l,n\},
\end{equation*}%
for any $z\in \left\{ z\in \left( K^{\times }\right) ^{n}\mid f_{j,\tau
}(z)=0,\text{ }j\in I\right\} $. Then by (\ref{con1}), 
\begin{equation*}
rank_{K}\left[ \frac{\partial f_{j,a}}{\partial x_{i}}\left( z\right) \right]
_{\substack{ j\in I  \\ 1\leq i\leq n}}=rank_{K}\left[ \frac{\partial
f_{j,\tau }}{\partial x_{i}}\left( z\right) \right] _{\substack{ j\in I  \\ %
1\leq i\leq n}}=\min \{l,n\},
\end{equation*}%
for any $z$ in 
\begin{equation*}
\left\{ z\in \left( K^{\times }\right) ^{n}\mid f_{j,a}(z)=0,\text{ }j\in
I\right\} \supseteq \left\{ z\in \left( K^{\times }\right) ^{n}\mid
f_{j,a}(z)=0,\text{ }j=1,\ldots ,l\right\} .
\end{equation*}%
Therefore, $\boldsymbol{f}$ is non-degenerate in the sense of Khovanskii.
\end{proof}

\begin{example}
\label{example4}Let $\boldsymbol{f}(x,y)=\left(
x^{2}-y^{2},x^{n},y^{m}\right) $, with $n,m\geq 3$. Then $\boldsymbol{f}$ is
not strongly non-degenerate at the origin with respect to $\Gamma (%
\boldsymbol{f})$. Indeed, $\Gamma (\boldsymbol{f})$ has only one compact
facet, $\tau $, that is the straight segment from $\left( 0,2\right) $\ to $%
\left( 2,0\right) $. Then 
\begin{equation*}
\boldsymbol{f}_{\tau }\left( x,y\right) \mathbf{=}\left(
x^{2}-y^{2},0,0\right) \text{, and \ }rank_{K}\left[ 
\begin{array}{cc}
2z_{1} & -2z_{2} \\ 
0 & 0 \\ 
0 & 0%
\end{array}%
\right] =1\neq \min \left\{ 2,3\right\} ,
\end{equation*}%
for every $\left( z_{1},z_{2}\right) \in \left\{ \left( z_{1},z_{2}\right)
\in \left( K^{\times }\right) ^{2}\mid z_{1}^{2}-z_{2}^{2}=0\right\} $, and
therefore $\boldsymbol{f}$ is not strongly non-degenerate with respect to $%
\Gamma (\boldsymbol{f})$. On the other hand, $\boldsymbol{f}$ is
non-degenerate in the sense of Khovanskii.
\end{example}

\subsection{Newton polyhedra and log-principalizations}

\begin{proposition}
\label{proposition1}Let $\boldsymbol{f}=(f_{1},\ldots ,f_{l}):U(\subseteq
K^{n})\longrightarrow K^{l}$ be a polynomial mapping (or more generally, an
analytic mapping defined \ on $U$) strongly \textit{non-degenerate at the
origin with respect to }$\Gamma \left( \boldsymbol{f}\right) $. Let $\ 
\mathcal{F}_{\boldsymbol{f}}$ be a simple fan subordinated to $\Gamma \left( 
\boldsymbol{f}\right) $. Let $Y_{K}$ be the toric manifold corresponding to $%
\mathcal{F}_{\boldsymbol{f}}$, and let \ 
\begin{equation*}
\sigma _{0}:Y_{K}\longrightarrow U
\end{equation*}%
be the restriction of the corresponding toric map to the inverse image of $U$%
. Denote by $Z$ the set of common zeroes of $\mathcal{I}_{\boldsymbol{f}%
}=(f_{1},\ldots ,f_{l})$ in $U\cap \left( K^{\times }\right) ^{n}$. When $U$
is taken small enough, either $Z=\varnothing $ or it is a submanifold of
codimension $l$. In this last case we have $l<n$ and we denote the closure
of $Z$ in $U$ and $Y_{K}$ by $Z_{U}$ and $Z_{Y}$, respectively.

\noindent (1) If $Z=\varnothing $ (or if $l=1$), the ideal $\sigma
_{0}^{\ast }\left( \mathcal{I}_{\boldsymbol{f}}\right) $ is principal (and
monomial)\ in a sufficiently small neighborhood \ of $\sigma
_{0}^{-1}\left\{ 0\right\} $.

\noindent (2) If $Z\neq \varnothing $, we have that $Z_{Y}$ is a closed
submanifold of $Y_{K}$, having normal crossings with the exceptional divisor
of $\sigma _{0}$. Let $\sigma _{1}:X_{K}\longrightarrow Y_{K}$ \ be \ the
blowing-up of $Y_{K}$\ with center $Z_{Y}$, and let $\sigma =\sigma
_{0}\circ \sigma _{1}:X_{K}\longrightarrow U$.\ Then the ideal $\sigma
^{\ast }\left( \mathcal{I}_{\boldsymbol{f}}\right) $ is principal\ (and
monomial) in a sufficiently small neighborhood \ of $\sigma ^{-1}\left\{
0\right\} $.
\end{proposition}

\begin{proof}
We first recall the construction of $\left( Y_{K},\sigma _{0}\right) $ from
a simple fan $\mathcal{F}_{\boldsymbol{f}}$ subordinated to $\Gamma \left( 
\boldsymbol{f}\right) $ (see e.g. \cite{AVG}). Let $\Delta _{\tau }$ be an $%
n-$dimensional \ simple cone in $\mathcal{F}_{\boldsymbol{f}}$ such that $%
F(a)=\tau $\ for \ any $a\in \Delta _{\tau }$. Then the face $\tau $ of $%
\Gamma \left( \boldsymbol{f}\right) $ is necessarily a point. Let $%
a_{1},\ldots ,a_{n}$ be the generators of $\Delta _{\tau }$. Then in the
chart of $Y_{K}$ corresponding to $\Delta _{\tau }$, the map $\sigma _{0}$
has the form

\begin{equation}
\begin{array}{cccc}
\sigma _{0}: & K^{n} & \longrightarrow & U \\ 
& y & \longrightarrow & x,%
\end{array}
\label{prin0}
\end{equation}%
where \ $x_{i}=\dprod\nolimits_{j}y_{j}^{a_{i,j}}$, with $\left[ a_{i,j}%
\right] =\left[ a_{1},\ldots ,a_{n}\right] $. Denote this chart by $V_{\tau
} $. We slightly abuse notation here : since $\sigma _{0}$ only maps to $U$
instead of to the whole of $K^{n}$, at some charts it will not be defined
everywhere on $K^{n}$. If $f_{i}\left( x\right) =\sum_{m}c_{m,i}x^{m}$ for $%
i=1,\ldots ,l$, then 
\begin{equation*}
\left( f_{i}\circ \sigma _{0}\right) \left( y\right)
=\sum_{m}c_{m,i}\dprod\limits_{j=1}^{n}y_{j}^{\left\langle
m,a_{j}\right\rangle }\text{ for \ }i=1,\ldots ,l.
\end{equation*}%
If supp$(f_{i})\cap \tau \neq \varnothing $, then the minimum of all $%
\left\langle m,a_{j}\right\rangle $ is attained at $\tau $, and then

\begin{equation}
\left( f_{i}\circ \sigma _{0}\right) \left( y\right) =\left(
\dprod\limits_{j=1}^{n}y_{j}^{d\left( a_{j}\right) }\right) \widetilde{f}%
_{i}\left( y\right) ,\text{ with }\widetilde{f}_{i}\left( 0\right) \neq 0
\label{prin2}
\end{equation}%
(cf. \cite[page 201, Lemma 8]{AVG}). If supp$(f_{i})\cap \tau =\varnothing $%
, 
\begin{equation}
\left( f_{i}\circ \sigma _{0}\right) \left( y\right) =\left(
\dprod\limits_{j=1}^{n}y_{j}^{d\left( a_{j}\right) }\right) \widetilde{f}%
_{i}\left( y\right) ,\text{ with }\widetilde{f}_{i}\left( 0\right) =0.
\label{prin3}
\end{equation}%
Then, from (\ref{prin2}) and (\ref{prin3}), we have in a neighborhood of the
origin of $V_{\tau }$ that $\sigma _{0}^{\ast }\left( \mathcal{I}_{%
\boldsymbol{f}}\right) $ is generated by $\tprod\nolimits_{j=1}^{n}y_{j}^{d%
\left( a_{j}\right) }$.

Now let us consider on $V_{\tau }$ the points on $\sigma _{0}^{-1}\left(
0\right) $, different from the origin of $V_{\tau }$. We will study
simultaneously points with exactly $r$ zero coordinates (where $1\leq r\leq
n-1$); after permuting indices, we may assume that the first $r$ coordinates
are zero.

Let $\tau ^{\prime }$ be the first meet locus of the cone $\Delta _{\tau
^{\prime }}$ spanned by $a_{1},\ldots ,a_{r}$; it is a compact face of $%
\Gamma \left( \boldsymbol{f}\right) $ (cf. \cite[page 201, Lemma 8]{AVG}).
We can write $\left( f_{i}\circ \sigma _{0}\right) \left( y\right) $ as

\begin{equation}
\left( f_{i}\circ \sigma _{0}\right) \left( y\right) =\left(
\dprod\limits_{j=1}^{r}y_{j}^{d\left( a_{j}\right) }\right) \left( 
\widetilde{f}_{i}\left( y_{r+1},\ldots ,y_{n}\right) +O_{i}(y_{1},\ldots
,y_{n})\right) ,  \label{prin4}
\end{equation}%
where the $\widetilde{f}_{i}$\ $\ $are polynomials in $y_{r+1},\ldots ,y_{n}$%
, and the $O_{i}(y_{1},\ldots ,y_{n})$ are analytic functions in $%
y_{1},\ldots ,y_{n}$ but belonging to the ideal generated by $y_{1},\ldots
,y_{r}$. Here the $\widetilde{f}_{i}$ are identically zero if and only if $%
\text{supp}(f_{i})\cap \tau ^{\prime }=\varnothing $. Furthermore,

\begin{equation}
(f_{i,\tau ^{\prime }}\circ \sigma _{0})\left( y\right) =\left(
\dprod\limits_{j=1}^{r}y_{j}^{d\left( a_{j}\right) }\right) \widetilde{f}%
_{i}\left( y_{r+1},\ldots ,y_{n}\right) .  \label{prin5a}
\end{equation}

We investigate the $\left( f_{i}\circ \sigma _{0}\right) \left( y\right) $
for $p=\left( 0,\dots ,0,p_{r+1},\dots ,p_{n}\right) $ with 
\begin{equation*}
\widetilde{p}=(p_{r+1},\dots ,p_{n})\in (K^{\times })^{n-r}.
\end{equation*}

We have to study two cases. The first case occurs when there exists an index 
$i$ such that $\widetilde{f}_{i}(\widetilde{p})\neq 0$. In this case as
before $\sigma _{0}^{\ast }\left( \mathcal{I}_{\boldsymbol{f}}\right) $ is
generated by $\tprod\nolimits_{j=1}^{r}y_{j}^{d\left( a_{j}\right) }$ in a
neighborhood of $p$.

The second case occurs when $\widetilde{f}_{i}(\widetilde{p})=0$, for all $%
i=1,\dots ,l$. We recall that, by the non-degeneracy condition, $rank_{K}%
\left[ \frac{\partial f_{i,\tau ^{\prime }}}{\partial x_{j}}\left( x\right) %
\right] =\min \{l,n\}$ for $x\in (K^{\times })^{n}\cap \{f_{1,\tau ^{\prime
}}(x)=\dots =f_{l,\tau ^{\prime }}(x)=0\}$. Since $\sigma _{0}$ is an
isomorphism over $\left( K^{\times }\right) ^{n}$, then also $rank_{K}\left[ 
\frac{\partial f_{i,\tau ^{\prime }}\circ \sigma _{0}}{\partial y_{j}}\left(
y\right) \right] $ $=\min \{l,n\}$ for $y\in (K^{\times })^{n}\cap
\{f_{1,\tau ^{\prime }}(\sigma _{0}(y))=\dots =f_{l,\tau ^{\prime }}(\sigma
_{0}(y))=0\}$. Note that by (\ref{prin5a}) this condition on $y$ is
equivalent to $y\in (K^{\times })^{n}\cap \{\widetilde{f_{1}}(y)=\dots =%
\widetilde{f_{l}}(y)=0\}$ and that $\left[ \frac{\partial f_{i,\tau ^{\prime
}}\circ \sigma _{0}}{\partial y_{j}}\left( y\right) \right] $ for such $y$
is equal to

\begin{equation*}
\left( 
\begin{array}{cccccc}
0 & \dots & 0 & \left( {\prod\nolimits_{j=1}^{r}}y_{j}^{d\left( a_{j}\right)
}\right) \frac{\partial \widetilde{f}_{1}}{\partial y_{r+1}}(y) & \dots & 
\left( {\prod\nolimits_{j=1}^{r}}y_{j}^{d\left( a_{j}\right) }\right) \frac{%
\partial \widetilde{f}_{1}}{\partial y_{n}}(y) \\ 
\dots & \dots & \dots & \dots & \dots & \dots \\ 
0 & \dots & 0 & \left( {\prod\nolimits_{j=1}^{r}}y_{j}^{d\left( a_{j}\right)
}\right) \frac{\partial \widetilde{f}_{l}}{\partial y_{r+1}}(y) & \dots & 
\left( {\prod\nolimits_{j=1}^{r}}y_{j}^{d\left( a_{j}\right) }\right) \frac{%
\partial \widetilde{f}_{l}}{\partial y_{n}}(y)%
\end{array}%
\right) .
\end{equation*}%
Now this implies that for $\widetilde{y}=(y_{r+1},\dots ,y_{n})\in
(K^{\times })^{n-r}\cap \{\widetilde{f_{1}}(\widetilde{y})=\dots =\widetilde{%
f_{l}}(\widetilde{y})=0\}$ the rank of the matrix

\begin{equation*}
\left( 
\begin{array}{ccc}
\frac{\partial \widetilde{f}_{1}}{\partial y_{r+1}}(\widetilde{y}) & \dots & 
\frac{\partial \widetilde{f}_{1}}{\partial y_{n}}(\widetilde{y}) \\ 
\dots & \dots & \dots \\ 
\frac{\partial \widetilde{f}_{l}}{\partial y_{r+1}}(\widetilde{y}) & \dots & 
\frac{\partial \widetilde{f}_{l}}{\partial y_{n}}(\widetilde{y})%
\end{array}%
\right)
\end{equation*}%
is equal to $\min \{l,n\}$. Then necessarily the rank is $l$, and we must
have that $l\leq n-r$.

So when $p$ above satisfies $\widetilde{f}_{i}(\widetilde{p})=0$ for $%
i=1,\dots ,l$, then necessarily \textit{all} $\widetilde{f}_{i}$ are nonzero
polynomials, $r\leq n-l$, and $rank_{K}\left[ \frac{\partial \widetilde{f}%
_{i}}{\partial y_{j}}\left( \widetilde{p}\right) \right] =l$. Now $\left[ 
\frac{\partial \widetilde{f}_{i}}{\partial y_{j}}\left( \widetilde{p}\right) %
\right] =\left[ \frac{\partial (\widetilde{f}_{i}+O_{i})}{\partial y_{j}}%
\left( p\right) \right] $ (cf. (\ref{prin4})). This last matrix having rank $%
l$ implies that we can choose new coordinates $y^{\prime }=(y_{1},\dots
,y_{r},y_{r+1}^{\prime },\dots ,y_{n}^{\prime })$ in a neighborhood $V_{p}$
of $p$ such that

\begin{equation}
\left( f_{i}\circ \sigma _{0}\right) \left( y^{\prime }\right) =\left(
\dprod\limits_{j=1}^{r}y_{j}^{d\left( a_{j}\right) }\right) y_{r+i}^{\prime }%
\text{ for }i=1,\dots ,l.  \label{prin6}
\end{equation}%
Since $\sigma _{0}$ is an isomorphism on $\left( K^{\times }\right) ^{n}$,
we have that $\{y_{r+1}^{\prime }=\dots =y_{r+l}^{\prime }=0\}$ is the
description in $V_{p}$ of $Z_{Y}\subset Y$. (The local description (\ref%
{prin6}) yields that $Z$ is a submanifold of $\left( K^{\times }\right) ^{n}$
of codimension $l$.) Clearly $Z_{Y}$ is a submanifold of $Y$ of codimension $%
l$, having normal crossings with the exceptional divisor of $\sigma _{0}$.

So, $\sigma _{1}$ being the blowing-up of $Y$ in $Z_{Y}$, we obtain by (\ref%
{prin6}) that $(\sigma _{0}\circ \sigma _{1})^{\ast }(\mathcal{I}_{%
\boldsymbol{f}})$ becomes principal.
\end{proof}

\begin{remark}
\label{global}If we replace in Proposition \ref{proposition1} the condition 
\textit{strongly non-degenerate at the origin with respect to }$\Gamma
\left( \boldsymbol{f}\right) $ by the condition \textit{strongly
non-degenerate with respect to }$\Gamma \left( \boldsymbol{f}\right) $, and $%
U$ by $K^{n}$, with a similar proof \ we obtain \ a \textit{global} version
of the proposition, that is, the conclusions (1) and (2) \ are valid without
the condition \textit{in a sufficiently small neighborhood}. In this case $%
Z_{Y}$ may have components that are disjoint with the exceptional divisor of 
$\sigma _{0}$.
\end{remark}

Given $\xi =\left( \xi _{1},\ldots ,\xi _{n}\right) \in \mathbb{N}%
^{n}\setminus \left\{ 0\right\} $, we put $\sigma \left( \xi \right) :=\xi
_{1}+\ldots +\xi _{n}$ and $d\left( \xi \right) =\min_{x\in \Gamma \left( 
\boldsymbol{f}\right) }\left\langle \xi ,x\right\rangle $ as before. We say
that $\xi $ is a primitive vector, if $\gcd \left( \xi _{1},\ldots ,\xi
_{n}\right) =1$. If $d\left( \xi \right) \neq 0$, we define 
\begin{equation*}
\mathcal{P}\left( \xi \right) =\left\{ -\frac{\sigma \left( \xi \right) }{%
d\left( \xi \right) }+\frac{2\pi \sqrt{-1}k}{d\left( \xi \right) \log q},%
\text{ \ }k\in \mathbb{Z}\right\} .
\end{equation*}

Let $\ \mathcal{F}_{\boldsymbol{f}}$ be a simple fan subordinated to $\Gamma
\left( \boldsymbol{f}\right) $. Then the set of generators of the cones in $%
\mathcal{F}_{\boldsymbol{f}}$, i.e. the skeleton of $\mathcal{F}_{%
\boldsymbol{f}}$, can be partitioned as $\Lambda _{\boldsymbol{f}}\cup 
\mathfrak{D}(\Gamma \left( \boldsymbol{f}\right) )$, where $\Lambda _{%
\boldsymbol{f}}$ is a finite set of primitive vectors, corresponding to the
extra rays, induced by the subdivision \ into simple cones.

The numerical data of the log-principalizations constructed in Proposition %
\ref{proposition1} and Remark \ref{global} can be computed directly from the
explicit expressions for the generators of $\sigma _{0}^{\ast }\left(
I_{f}\right) $, $\sigma ^{\ast }\left( I_{f}\right) $, and Lemma 8 in \cite[%
page 201]{AVG}. Then Theorem \ref{Theorem1}\ \ yields that the poles of $%
Z_{\Phi }\left( s,\boldsymbol{f}\right) $ \ belong to the set 
\begin{equation}
\tbigcup\limits_{\xi \in \Lambda _{\boldsymbol{f}}}\mathcal{P}\left( \xi
\right) \cup \tbigcup_{\xi \in \mathfrak{D}(\Gamma \left( \boldsymbol{f}%
\right) )}\mathcal{P}\left( \xi \right) \cup \left\{ -l+\frac{2\pi \sqrt{-1}k%
}{\log q}\text{, }k\in \mathbb{Z}\right\} ,  \label{list}
\end{equation}%
where the last set may be discarded if $l\geq n$.

This provides a generalization to the case \ $l\geq 1$ of a well-known
result that describes the poles of the local zeta function associated to a
non-degenerate polynomial in terms of the corresponding Newton polyhedron 
\cite{L-M}, \cite{D1a}, \cite{D-H}, \cite{Z1}. This result was originally
established by Varchenko \cite{Var}\ for local zeta functions over $\mathbb{R%
}$. As in the case $l=1$, the list (\ref{list}) is too big. More precisely,
the set $\cup _{\xi \in \Lambda _{\boldsymbol{f}}}\mathcal{P}\left( \xi
\right) $ is not necessary. This fact is established by analogous arguments
\ as in \cite{D1a} where the case $l=1$ is studied.

\begin{theorem}
\label{Theorem2a}(1) Let $\boldsymbol{f}=(f_{1},\ldots
,f_{l}):U\longrightarrow K^{l}$ be an analytic\textit{\ mapping} strongly 
\textit{non-degenerate \ at the origin with respect to }$\Gamma \left( 
\boldsymbol{f}\right) $\textit{. }If $U$ is a sufficiently small
neighborhood of the origin, and $\Phi $ is a \ Schwartz-Bruhat \ function \
whose support \ is contained in \ $U$, then the poles of $Z_{\Phi }\left( s,%
\boldsymbol{f}\right) $ belong to the set $\cup _{\xi \in \mathfrak{D}%
(\Gamma \left( \boldsymbol{f}\right) )}\mathcal{P}\left( \xi \right) \cup
\left\{ -l+\frac{2\pi \sqrt{-1}k}{\log q}\text{, }k\in \mathbb{Z}\right\} $,
where the last set may be discarded \ if $l\geq n$.

\noindent (2) If $\boldsymbol{f}$ $:K^{n}\longrightarrow K^{l}$\textit{\ is
a strongly non-degenerate polynomial mapping with respect to }$\Gamma \left( 
\boldsymbol{f}\right) $, then the poles of $Z\left( s,\boldsymbol{f}\right) $
belong to the set 
\begin{equation*}
\cup _{\xi \in \mathfrak{D}(\Gamma \left( \boldsymbol{f}\right) )}\mathcal{P}%
\left( \xi \right) \cup \left\{ -l+\frac{2\pi \sqrt{-1}k}{\log q}\text{, }%
k\in \mathbb{Z}\right\} .
\end{equation*}
\end{theorem}

The above result can be restated in a geometric form as follows. If $s$ is a
pole of $Z_{\Phi }\left( s,\boldsymbol{f}\right) $, then $\func{Re}(s)$ is $%
-l$, or $\func{Re}(s)$ is \ of the form $-1/t_{0}$, where $\left(
t_{0},\ldots ,t_{0}\right) $ \ is the intersection \ point \ of the diagonal 
$\{\left( t,\ldots ,t\right) \in \mathbb{R}^{n}\}$ with the supporting
hyperplane of a facet of $\Gamma \left( \boldsymbol{f}\right) $.

By using Theorems \ref{theorem1c} and \ref{Theorem2a} we obtain the
following corollary.

\begin{corollary}
\label{coro2b}(1) Let $U$ be a sufficiently small neighborhood of the
origin, and let $\boldsymbol{f}=(f_{1},\ldots ,f_{l}):U\longrightarrow K^{l}$
be an analytic\textit{\ mapping} strongly \textit{non-degenerate \ at the
origin with respect to }$\Gamma \left( \boldsymbol{f}\right) $\textit{. }Let 
$\left( t_{\boldsymbol{f}},\ldots ,t_{\boldsymbol{f}}\right) \in \mathbb{Q}%
^{n}$ be the intersection \ point \ of the diagonal $\{\left( t,\ldots
,t\right) \in \mathbb{R}^{n}\}$ with the boundary of $\Gamma \left( \mathbf{f%
}\right) $. If $t_{\boldsymbol{f}}\geq 1/l$, then \ $-1/t_{\boldsymbol{f}}$
is the largest real part of a pole of $Z_{U}\left( s,\boldsymbol{f}\right) $.

\noindent (2) Let $\boldsymbol{f}$ $:K^{n}\longrightarrow K^{l}$\textit{\ be
a strongly non-degenerate polynomial mapping with respect to }$\Gamma \left( 
\boldsymbol{f}\right) $. If $t_{\boldsymbol{f}}\geq 1/l$, then \ $-1/t_{%
\boldsymbol{f}}$ is the largest real part of a pole of $Z\left( s,%
\boldsymbol{f}\right) $.
\end{corollary}

The largest real part of the poles of $Z\left( s,\boldsymbol{f}\right) $, $%
l=1$, when $\boldsymbol{f}$ is non-degenerate with respect to its Newton
polyhedron $\Gamma \left( \boldsymbol{f}\right) $ and $t_{\boldsymbol{f}}>1$
follows from observations made by Varchenko in \cite{Var} and was originally
noted in the $p-$adic case in \cite{L-M}. The case $t_{\boldsymbol{f}}=1$ is
treated in \cite{D-H}. The case of $t_{\boldsymbol{f}}<1$ is more difficult
and is established in \cite{D-H} with some additional conditions on $\Gamma
\left( \boldsymbol{f}\right) $ by using a difficult result on exponential
sums. In \cite{Z1} the second author established the case $t_{\boldsymbol{f}%
}\geq 1$ when $\boldsymbol{f}$ \ is a non-degenerate polynomial with
coefficients in a non-archimedean local field of arbitrary characteristic.

\section{Explicit formulas and Newton polyhedra}

In \cite[Theorem 4.2]{D-H} Denef and Hoornaert gave an explicit formula for $%
Z(s,\boldsymbol{f})$, $l=1$, associated to a polynomial $\boldsymbol{f}$%
\textbf{\ }in several variables over the $p-$adic numbers, when $\boldsymbol{%
f}$ is sufficiently non-degenerate with respect to its Newton polyhedron $%
\Gamma \left( \boldsymbol{f}\right) $. This explicit formula can be
generalized to the case $l\geq 1$ \ by using the \ condition of
non-degeneracy for polynomial mappings introduced in this paper.

Let as before $K$ \ be \ a $p-$adic field \ with valuation ring $R_{K}$,
maximal ideal $P_{K}$ and residue field $\overline{K}=\mathbb{F}_{q}$. For
any polynomial $g$ over $R_{K}$ we denote by $\overline{g}$ \ the polynomial
over $\overline{K}$ obtained by reducing each coefficient of $g$ modulo $%
P_{K}$.

\begin{definition}
\label{def2} Let $f_{i}\in R_{K}\left[ x\right] $, $x=\left( x_{1},\ldots
,x_{n}\right) $, satisfying $f_{i}\left( 0\right) =0$ for $i=1,\ldots ,l$.
The mapping $\boldsymbol{f}=(f_{1},\ldots ,f_{l}):K^{n}\longrightarrow K^{l}$
\ is called \textit{strongly} \textit{non-degenerate over }$\overline{K}$%
\textit{\ with respect to }$\Gamma (\boldsymbol{f})$, if for \textit{\ }any
face $\tau $ of $\Gamma (\boldsymbol{f})$, including $\Gamma (\boldsymbol{f}%
) $ itself, we have that $rank_{K}\left[ \frac{\partial \overline{%
f_{i,_{\tau }}}}{\partial x_{j}}\left( \overline{z}\right) \right] =\min
\left\{ l,n\right\} $, for any $\overline{z}\in \left( \overline{K}^{\times
}\right) ^{n}$ satisfying $\overline{f_{1,\tau }}(\overline{z})=\ldots =%
\overline{f_{l,\tau }}(\overline{z})=0$. Analogously we call $\boldsymbol{f}$
\textit{strongly non-degenerate at the origin over }$\overline{K}$\textit{\
with respect to }$\Gamma (\boldsymbol{f})$, if the same condition is
satisfied but only for the \textit{compact} faces $\tau $ \ of $\Gamma (%
\boldsymbol{f}) $.
\end{definition}

\begin{theorem}
\label{Theorem2} (1) Let $\boldsymbol{f}=(f_{1},\ldots
,f_{l}):K^{n}\rightarrow K^{l}$ be a strongly non-degenerate polynomial
mapping over $\overline{K}$. Denote \ for each \ face $\tau $ of $\Gamma (%
\boldsymbol{f})$, including $\Gamma (\boldsymbol{f})$ itself, 
\begin{equation*}
\overline{D}_{\tau }:=\left\{ \overline{x}\in \left( \overline{K}^{\times
}\right) ^{n}\mid \overline{f_{1,\tau }}\left( \overline{x}\right) =\ldots =%
\overline{f_{l,\tau }}\left( \overline{x}\right) =0\right\} .
\end{equation*}%
Fix a rational simplicial polyhedral subdivision $\left\{ \Delta _{\tau
,i}\right\} $, with $\tau $ a proper face,\ subordinated to $\Gamma \left( 
\boldsymbol{f}\right) $ as in (\ref{simplisubv}). Denote by $a_{j}$, $%
j=1,\ldots $, $r_{\Delta _{\tau ,i}}$, the generators of the \ cone $\Delta
_{\tau ,i}$. Then 
\begin{equation*}
Z(s,f)=L_{\Gamma \left( \boldsymbol{f}\right) }\left( q^{-s}\right)
+\dsum\limits_{\tau \neq \Gamma \left( \boldsymbol{f}\right) }L_{\tau
}\left( q^{-s}\right) \left( \tsum\limits_{i}S_{\tau ,i}\left( q^{-s}\right)
\right) .
\end{equation*}%
Here 
\begin{equation*}
L_{\tau }\left( q^{-s}\right) =q^{-n}\left( (q-1)^{n}-\frac{\text{card}(%
\overline{D}_{\tau })\left( 1-q^{-s}\right) }{1-q^{-\min \left\{ l,n\right\}
-s}}\right) ,
\end{equation*}%
for each face $\tau $ of $\Gamma \left( \boldsymbol{f}\right) $, including $%
\Gamma \left( \boldsymbol{f}\right) $, \ and%
\begin{equation*}
S_{\tau ,i}\left( q^{-s}\right) =\frac{\left( \tsum\limits_{h}q^{\sigma
\left( h\right) +d\left( h\right) s}\right) q^{^{-\sum_{j=1}^{r_{\Delta
_{\tau ,i}}}\left( \sigma \left( a_{j}\right) +d\left( a_{j}\right) s\right)
}}}{\dprod\limits_{j=1}^{r_{\Delta _{\tau ,i}}}\left( 1-q^{-\sigma \left(
a_{j}\right) -d\left( a_{j}\right) s}\right) },
\end{equation*}%
where $h$ runs through \ the elements of the set 
\begin{equation*}
\mathbb{Z}^{n}\cap \left\{ \tsum\limits_{j=1}^{r_{\Delta _{\tau ,i}}}\lambda
_{j}a_{j}\mid 0\leq \lambda _{j}<1\text{ for }j=1,\ldots ,r_{\Delta _{\tau
,i}}\right\} .
\end{equation*}%
(2) With the same \ notations and only assuming that \ $\boldsymbol{f}$ is
strongly non-degenerate at the origin over $\overline{K}$ we have%
\begin{equation*}
Z_{0}\left( s,\boldsymbol{f}\right) =\dsum\limits_{\tau \text{ compact}%
}L_{\tau }\left( q^{-s}\right) \left( \tsum\limits_{i}S_{\tau ,i}\left(
q^{-s}\right) \right) .
\end{equation*}
\end{theorem}

The proof of the above result is analogous to the case $l=1$ treated in \cite%
[Theorem 4.2]{D-H}.

By using a simple polyhedral subdivision one obtains a slightly less
complicated explicit formula \ in which all the terms $\tsum\nolimits_{h}q^{%
\sigma \left( h\right) +d\left( h\right) s}$ are identically $1$. But then
in general we have to introduce new rays which give rise to superfluous
candidate poles.

\begin{example}
\label{ex5}Let $\boldsymbol{f}=\left( x^{3}-xy,y\right) $ as in Example \ref%
{Exam4}. It is strongly non-degenerate over $\overline{K}$ with respect to $%
\Gamma (\boldsymbol{f})$. We shall compute $Z(s,\boldsymbol{f})$ using
Theorem \ref{Theorem2}\ and the obvious rational simplicial \ polyhedral
subdivision of $\mathbb{R}_{+}^{2}$. More precisely, set $a_{1}=\left(
0,1\right) $, $a_{2}=\left( 1,3\right) $, and $a_{3}=\left( 1,0\right) $; $%
\Delta _{i}=\left\{ a_{i}\lambda \mid \lambda >0\right\} $ for $i=1$, $2$, $3
$, and $\Delta _{i,i+1}=\left\{ \lambda a_{i}+\lambda ^{\prime }a_{i+1}\mid
\lambda ,\lambda ^{\prime }>0\right\} $, $i=1$, $2$. Then 
\begin{equation*}
\mathbb{R}_{+}^{2}=\left\{ 0\right\} \cup \Delta _{1}\cup \Delta _{1,2}\cup
\Delta _{2}\cup \Delta _{2,3}\cup \Delta _{3}.
\end{equation*}%
\ With the notation of Theorem \ref{Theorem2}\ \ one easily verifies that
all $\overline{D}_{\tau }=\varnothing $ and hence all $L_{\tau
}=q^{-2}\left( q-1\right) ^{2}$. Further 
\begin{equation*}
S_{\tau _{1}}=S_{\tau _{3}}=\frac{q^{-1}}{1-q^{-1}}\text{, \ }S_{\tau _{2}}=%
\frac{q^{-4-3s}}{1-q^{-4-3s}}\text{,}
\end{equation*}%
\begin{equation*}
S_{\tau _{1,2}}=\frac{q^{-5-3s}}{\left( 1-q^{-1}\right) \left(
1-q^{-4-3s}\right) }\text{, \ }S_{\tau _{2,3}}=\frac{\left(
1+q^{2+s}+q^{3+2s}\right) q^{-5-3s}}{\left( 1-q^{-1}\right) \left(
1-q^{-4-3s}\right) }.
\end{equation*}%
Therefore%
\begin{equation*}
Z(s,\boldsymbol{f})=q^{-2}\left( q-1\right) \frac{\left(
q+1+q^{-1-s}+q^{-2-2s}\right) }{1-q^{-4-3s}}.
\end{equation*}%
If we would use the natural \textit{simple} polyhedral subdivision of the
one above, introducing two new rays generated by $\left( 1,1\right) $ and $%
\left( 1,2\right) $, we would introduce the same \ superfluous (real)
candidate poles $-2$ and $-\frac{3}{2}$ as in Example \ref{Exam4}. This is
reasonable \ because the log-principalization of Proposition \ref%
{proposition1} associated to this simple fan is in fact the same as the one
constructed in Example \ref{Exam4}.
\end{example}

\begin{example}
Let $\boldsymbol{f}=\left( y^{2}-x^{3},y^{2}-z^{2}\right) $ as in Example %
\ref{Example5}. When $char\left( \overline{K}\right) \neq 2$, it is strongly
non-degenerate at the origin over $\overline{K}$ with respect to $\Gamma (%
\boldsymbol{f})$. The Newton polyhedron $\Gamma (\boldsymbol{f})$ has seven
compact faces.\ The polyhedral subdivision associated to \ it is already
simplicial, so in the formulation of Theorem \ref{Theorem2}\ (2) we need to
sum over seven cones: the ray through $a=\left( 2,3,3\right) $, the three $2-
$dimensional cones with \ $a$ in their boundaries, and the three $3-$%
dimensional cones. We note that all the $\overline{D}_{\tau }=\varnothing $,
except \ when $\tau $\ is the unique compact facet, in this case card$\left( 
\overline{D}_{\tau }\right) =2\left( q-1\right) $. Concerning the $S_{\tau
}\left( q^{-s}\right) $ we just mention that \ the expression $%
\tsum\nolimits_{h}q^{\sigma \left( h\right) +d\left( h\right) s}$ is three
times equal to $1$, three times equal to $1+q^{3+2s}+q^{6+4s}$, and once to $%
1+q^{5+3s}$. One can verify \ that the formula \ in Theorem \ref{Theorem2}\
\ yields the same expression for $Z_{0}(s,\boldsymbol{f})$ as in Example \ref%
{Example5}. Note that $-8/6$ and $-2$ are the only \ (real) candidate poles
given by Theorems \ref{Theorem2a}\ or \ref{Theorem2}.
\end{example}

\begin{remark}
With the obvious analogous definitions for strongly non-degeneracy over $%
\mathbb{C}$, we have \ the following. Suppose \ that $f_{1},\ldots f_{l}$
are polynomials in $n$ variables with coefficients in a number field $F$ $%
\left( \subseteq \mathbb{C}\right) $. Then we can consider $\boldsymbol{f}%
=(f_{1},\ldots ,f_{l})$ as a map $K^{n}\rightarrow K^{l}$ for any
non-archimedean completion $K$ of $F$. If $\boldsymbol{f}$ is strongly
non-degenerate at the origin over $\mathbb{C}$\ with respect to $\Gamma
\left( \boldsymbol{f}\right) $, then $\boldsymbol{f}$ is strongly
non-degenerate over $\overline{K}$\ with respect to $\Gamma \left( 
\boldsymbol{f}\right) $\ for almost all the completions \ $K$ of $F$. (And
analogously for non-degeneracy at the origin.) \ This \ fact follows by
applying the Weak Nullstellensatz.
\end{remark}

\begin{remark}
By using our notion of non-degeneracy with respect to a Newton polyhedron it
is also possible to give lists of candidate poles and explicit formulae for
the motivic and topological zeta functions \ introduced in \ref{motivic},
associated to a polynomial ideal. These explicit formulas are reasonably
straightforward \ generalizations of those in \cite{A-C-L-M} and \ \cite[Th%
\'{e}or\`{e}me 5.3 (i)]{D-L3}. For the topological zeta function one
requires here strongly non-degeneracy with respect to all the faces of the
\textquotedblleft global\textquotedblright\ Newton \ polyhedron as in \cite[%
(5.1)]{D-L3}.
\end{remark}

\end{document}